\documentclass[12pt]{amsart}
\usepackage[all, knot]{xy}
\usepackage{epsfig}
\usepackage{xspace}
\usepackage{layout}
\usepackage{latexsym}
\usepackage{booktabs, amsmath}
\usepackage{amsthm}
\usepackage{amssymb}
\usepackage{amsfonts}
\usepackage{amsxtra}     
\usepackage{color}
\usepackage{verbatim}
\usepackage{longtable}
\usepackage{url}
\usepackage{tikz}
\usetikzlibrary{matrix,arrows,decorations.pathmorphing}

\newcommand{\co}{\colon\thinspace}

\newcounter{commentcounter}

\begin{document}

\newtheorem{thm}{Theorem}[section]
\newtheorem{conj}[thm]{Conjecture}
\newtheorem{lem}[thm]{Lemma}
\newtheorem{cor}[thm]{Corollary}
\newtheorem{prop}[thm]{Proposition}
\newtheorem{rem}[thm]{Remark}

\theoremstyle{definition}
\newtheorem{defn}[thm]{Definition}
\newtheorem{examp}[thm]{Example}
\newtheorem{construction}[thm]{Construction}
\newtheorem{notation}[thm]{Notation}
\newtheorem{rmk}[thm]{Remark}

\theoremstyle{remark}

\makeatletter
\renewcommand{\maketag@@@}[1]{\hbox{\m@th\normalsize\normalfont#1}}%
\makeatother

\renewcommand{\labelenumi}{(\roman{enumi})}
\renewcommand{\labelenumii}{(\alph{enumii})}

\renewcommand{\theenumi}{(\roman{enumi})}
\renewcommand{\theenumii}{(\alph{enumii})}

\def\square{\hfill ${\vcenter{\vbox{\hrule height.4pt \hbox{\vrule width.4pt
height7pt \kern7pt \vrule width.4pt} \hrule height.4pt}}}$}

\newenvironment{pf}{{\it Proof:}\quad}{\square \vskip 12pt}
\date{\today}

\title[GIS Morphisms and $PR$-Groups]{Morphisms of Generalized Interval Systems and $PR$-Groups}
\author[Fiore, Noll, Satyendra]{Thomas M. Fiore, Thomas Noll, and Ramon Satyendra}
\address{Thomas M. Fiore \\ Department of Mathematics and Statistics\\
University of Michigan-Dearborn \\ 4901 Evergreen Road \\ Dearborn,
MI 48128 \\ U.S.A.} \email{tmfiore@umich.edu}
\urladdr{http://www-personal.umd.umich.edu/~tmfiore/}
\address{Thomas Noll \\ Escola Superior de M\'{u}sica de Catalunya \\
Departament de Teoria, Composici\'{o} i Direcci\'{o} \\
C. Padilla, 155 - Edifici L'Auditori \\
08013 Barcelona, Spain }
\email{noll@cs.tu-berlin.de}
\urladdr{http://user.cs.tu-berlin.de/~noll/}
\address{Ramon Satyendra \\ School of Music, Theatre and Dance \\
University of Michigan \\ 1100 Baits Drive \\ Ann Arbor,
MI 48109-2085 \\ U.S.A.} \email{ramsat@umich.edu}
\urladdr{http://ramonsatyendra.net/}


\begin{abstract}
We begin the development of a categorical perspective on the theory of generalized interval systems (GIS's). Morphisms of GIS's allow the analyst to move between multiple interval systems and connect transformational networks. We expand the analytical reach of the Sub Dual Group Theorem of Fiore--Noll \cite{fiorenoll2011} and the generalized contextual group of Fiore--Satyendra \cite{fioresatyendra2005} by combining them with a theory of GIS morphisms. Concrete examples include an analysis of Schoenberg, String Quartet in $D$ minor, op. 7, and simply transitive covers of the octatonic set. This work also lays the foundation for a transformational study of Lawvere--Tierney upgrades in the topos of triads of Noll \cite{nollToposOfTriads}.
\end{abstract}

\maketitle

\tableofcontents

\section{Introduction}

In this article, we introduce morphisms of generalized interval systems in an analysis of Schoenberg\footnote{Arnold
Sch\"onberg, Quartett f\"ur zwei Violinen, Viola und Violoncello, Opus 7,
Berlin: Verlag Dreililien, [1905].}, String Quartet in $D$ minor, op. 7. This categorical perspective on generalized interval systems allows us to expand the analytical reach of the Sub Dual Group Theorem of \cite{fiorenoll2011} and the generalized contextual group of \cite{fioresatyendra2005} to move between generalized interval systems, connect transformational networks, and circumvent the stringent condition of simple transitivity without giving up its virtues.

\begin{figure}[h]\caption{Schoenberg, String Quartet in $D$ minor, op. 7. Complete $RP$-chain of triads in measures 88-92} \label{fig:Figure5CompleteRPChain88-92}
\includegraphics[height=1.4in]{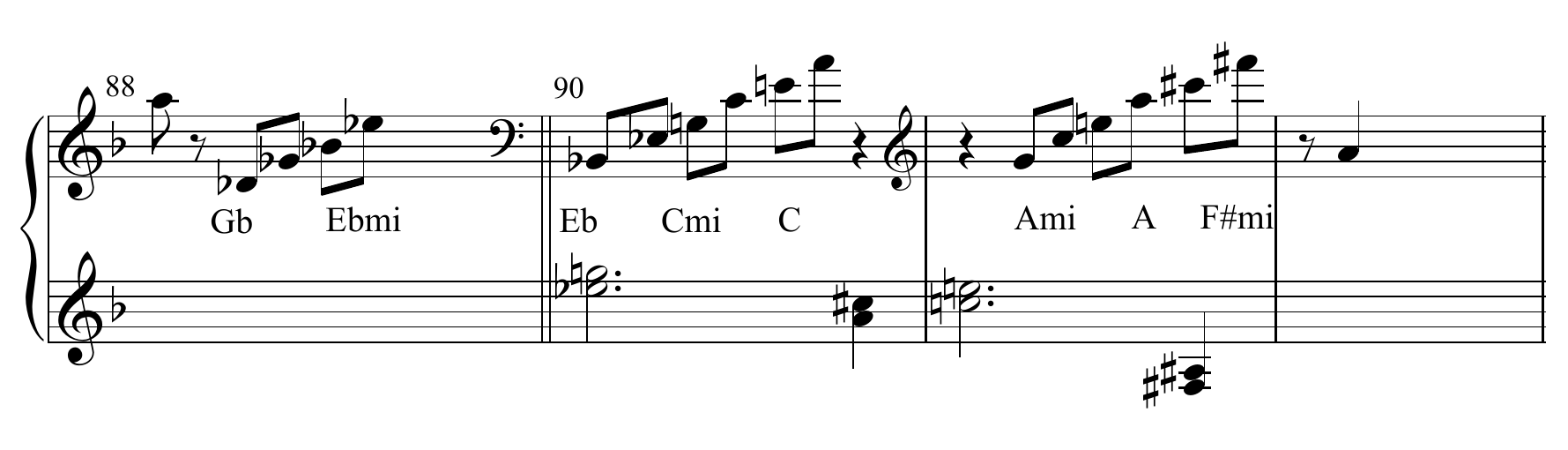}
\end{figure}

To motivate our development of GIS morphisms, consider the
triadic melody in Figure~\ref{fig:Figure5CompleteRPChain88-92}, which ends the first large formal section of the quartet. These consonant triads form a complete $RP$-chain, that is to say, they are an orbit of the subgroup $\langle R, P \rangle$ of the neo-Riemannian $PLR$-group, formed by iteratively applying $R$ then $P$. Furthermore, these consonant triads form a maximal cover of the octatonic: their union is the pitch-class set $\{C, C\sharp, E\flat, E, G\flat, G, A,B\flat\}$, and there are no other consonant triads in this pitch-class set besides these listed in Figure~\ref{fig:Figure5CompleteRPChain88-92}.  For cardinality reasons, the eight-element subgroup $\langle R, P \rangle$ acts simply transitively on this set of consonant triads. We thus have a {\it sub generalized interval system} of the full neo-Riemannian $PLR$-group action on all 24 consonant triads.

\begin{figure}\caption{Equivalence of the triadic cycle to a shark-jet cycle via the affine transformation ${}^77$} \label{fig:Figure6TransformationOfPRCycleNotes}
\includegraphics[height=1.75in]{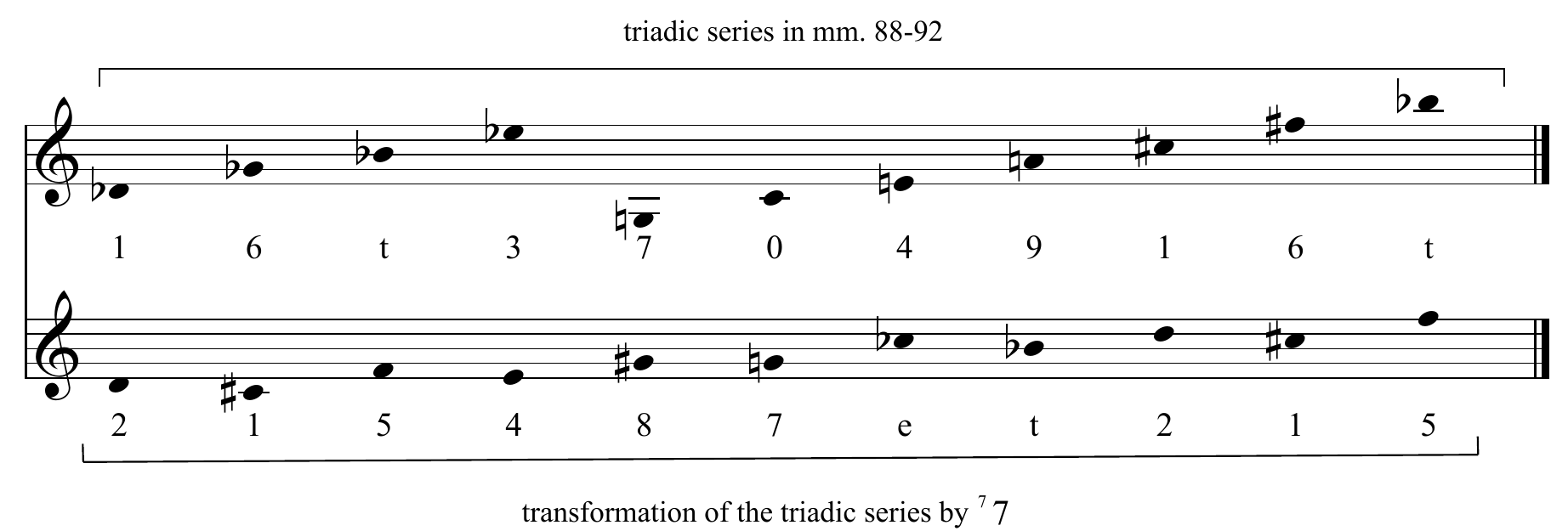}
\end{figure}

\begin{figure}\caption{Affine image of triadic melody from mm. 88-92 as a piece-wide narrative constructed from the opening motivic cell from measures 1-2} \label{fig:Figure7PiecewideNarrative}
\includegraphics[height=1.3in]{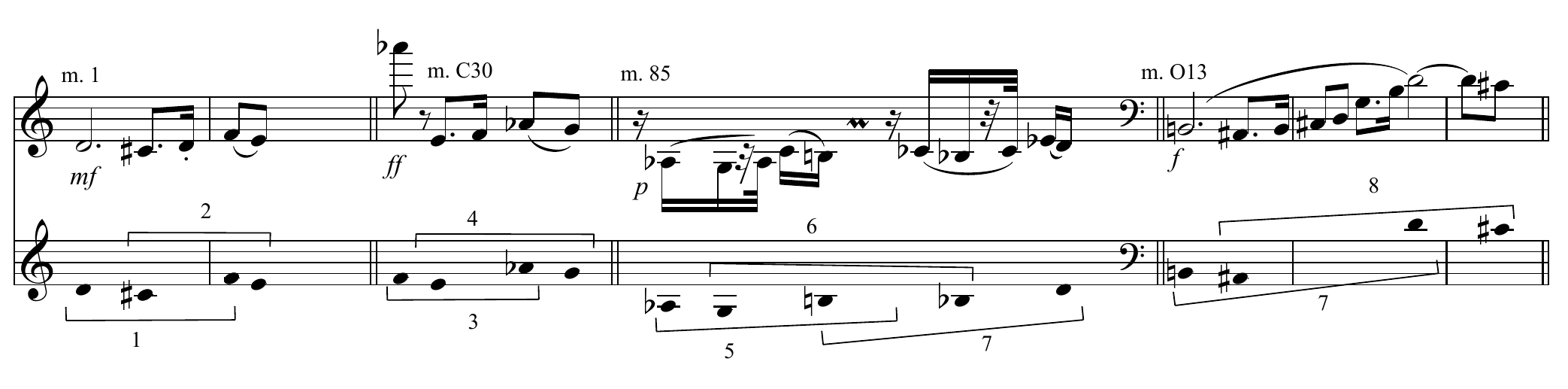}
\end{figure}

\begin{figure} \caption{Affine map ${}^77$ sends interlocking consonant triads to interlocking jet-shark trichords}
\label{fig:77_commutativity}
\begin{math}
\renewcommand{\labelstyle}{\textstyle}
\entrymodifiers={+[F-]}
\xymatrix@C=4pc@R=3.5pc{(1,6,10) \ar[r]^{(13)R} \ar[d]_{{}^77} &  (6,10,3) \ar[d] \ar[d]^{{}^77}  \\ (2,1,5) \ar[r]^{(13)R}  & (1,5,4)}
\end{math}
\end{figure}
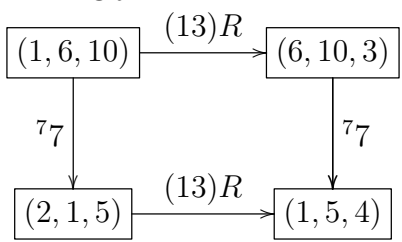

Next, consider the image of this triadic cycle under the affine map ${}^77 \colon x \mapsto 7x+7$ pictured in bottom staff of Figure~\ref{fig:Figure6TransformationOfPRCycleNotes}, and displayed as a piece-wide narrative in Figure~\ref{fig:Figure7PiecewideNarrative}.  This image consists of a maximal cover of the octatonic $\{C\sharp,D,E,F,G,G\sharp,B\flat,C\flat\}$ by 8 elements of the $T/I$-class of $(2,1,5)$ (which we call {\it jets} and {\it sharks} because of their role in Bernstein's {\it West Side Story}). Every three consecutive notes $s_i$, $s_{i+1}$, $s_{i+2}$ in the bottom staff form an element of the $T/I$-class of $(2,1,5)$, interlocked with the previous three consecutive notes  $s_{i-1}$, $s_{i}$, $s_{i+1}$; that is, we have another ``$PR$''-cycle in the bottom staff, but in the $(2,1,5)$-class. Note also that ${}^77$ maps interlocking chords to interlocking chords, so in the first four notes of the two staves for instance, we have the commutativity in Figure~\ref{fig:77_commutativity}. That is, the affine map ${}^77$ is a {\it morphism of sub generalized interval systems}.

We have now arrived at a list of desiderata for a relevant mathematical framework.\footnote{We do not claim that this list is universal or exhaustive.} Such a mathematical theory should:
\begin{enumerate}
\item
Associate a neo-Riemannian-type group to $3$-tuples $(x_1, x_2,x_3)$ of pitch-classes, e.g. for $(2,1,5)$ in the opening theme (see Figure~\ref{fig:Figure1OpeningTheme}) \smallskip \smallskip
\item
Have attendant theorems about duality between this neo-Riemannian type group and the $T/I$-group \smallskip \smallskip
\item
Show how one can move via morphisms between the associated neo-Riemannian type groups associated to {\it two} such $3$-tuples $(x_1, x_2,x_3)$ and $(y_1, y_2,y_3)$ when related by an affine map, e.g. relating via $^77$ the major $(1,6,10)$ to the jet $(2,1,5)$ trichords \smallskip\smallskip
\item
Show how to obtain substructures of neo-Riemannian type groups, e.g. for the octatonic and its maximal major/minor cover or maximal jet/shark cover\smallskip \smallskip
\item
Show how to move between these substructures, e.g. between the three octatonics\smallskip \smallskip
\item
Determine which subsets of the octatonic set generate maximal covers with simply transitive actions under neo-Riemannian type group actions. \smallskip \smallskip
\end{enumerate}

Neo-Riemannian type groups associated to pitch-class segments and substructures of dual groups have already appeared in the literature. Indeed these are the main topics of our ealier papers \cite{fioresatyendra2005} respectively \cite{fiorenoll2011}, as we recall in Section~\ref{sec:Sub_Dual_Groups_and_Contextual_Groups}. Morphisms of generalized interval systems, from both the intervallic and transformational points of view, are the main topic of the present paper, and are treated in Section~\ref{sec:Morphisms}.

In Section~\ref{sec:Octatonic_Example}, we use the foregoing theory to generate three dual pairs of groups acting on three maximal covers of the octatonic. These originate from the major $(0,4,7)$, the jet $(0,4,1)$, and the stride\footnote{Our use of the nickname {\it stride} was motivated by stride piano techniques, where this chord is often played in wide position in the left hand.} $(0,4,10)$.  Multiplication by 7 and multiplication by 10 are morphisms from the major system to the jet system and the stride system respectively. Interestingly, the latter morphism is neither injective nor surjective, and the stride system is {\it not} generated by $L$ and $R$-analogues. This detailed example paves the way for our analysis in Section~\ref{sec:Analysis_of_Schoenberg} of Schoenberg's Quartet in $D$ minor.

In Section~\ref{sec:Covers}, we show that any three-element subset of the octatonic gives rise to a simply transitive action by the set-wise $T/I$-stabilizer of the octatonic, though such a cover may not be a substructure of a $T/I$-contextual group pair because of symmetries in the three-element subset.

\noindent {\bf Acknowledgements.} The authors acknowledge the first author's Rackham Faculty Research Grant of the University of Michigan, which supported Thomas Noll's visit to Dearborn and Ann Arbor in September, 2011. The authors thank the referees for their helpful comments.

\section{Morphisms of Generalized Interval Systems} \label{sec:Morphisms}

We begin by developing morphisms of generalized interval systems first in an intervallic and then in a transformational context.
We then prove that the categories of generalized interval systems, simply transitive group actions, and groups (with affine morphisms) are all equivalent. We also characterize monic/epic/iso morphisms in the category of generalized interval systems.

In many works of music, there are several generalized interval systems upon which analysts draw, each with its own notion of interval and transposition. It is therefore desirable to have a codified way of moving between these systems, so that intervals correspond to intervals and transpositions to transpositions. For instance, one GIS may be contained in another, a situation we discuss in more detail in Section~\ref{sec:Sub_Dual_Groups_and_Contextual_Groups} from the transformational point of view. Another situation in which morphisms are useful is with products: each product of generalized interval systems projects to its factors. Morphisms also arise in the context of quotients:
if a generalized interval system has a congruence relation on its interval group in the sense of \cite[pages 34--35]{LewinGMIT}, then there is a morphism to the quotient GIS.

More specifically, the concrete application which motivates our development of GIS morphisms is our analysis in Section 5 of Schoenberg's Quartet in $D$ minor, as we have briefly sketched in the Introduction. In the transformational formulation of generalized interval systems, the analysis contains several morphisms. We have the inclusion of the simply transitive $PR$-group action on the 8 chords in Figure~\ref{fig:Figure5CompleteRPChain88-92} into the simply transitive  $PLR$-group action on the 24 major and minor triads. This inclusion morphism also has its (2,1,5)-counterpart, beginning with the 8 trichords in the bottom staff of Figure~\ref{fig:Figure6TransformationOfPRCycleNotes}.  More interesting is the affine map ${}^77$, which is part of a morphism from the
simply transitive $PR$-group action on the chords in Figure~\ref{fig:Figure5CompleteRPChain88-92} to another ``$PR$''-group action on 8 chords in the $(2,1,5)$-class, as pictured in Figure~\ref{fig:Figure6TransformationOfPRCycleNotes},  the piecewide narrative in Figure~\ref{fig:Figure7PiecewideNarrative}, and the example of the morphism property in Figure~\ref{fig:77_commutativity}. This morphism is even an isomorphism, as ${}^77$ is bijective (see Proposition~\ref{prop:monic/epic_morphisms_for_intervallic}).

Recall that a {\it generalized interval system} $(S,\textit{IVLS},\text{int})$, as invented by David Lewin in \cite{LewinGMIT}, consists of a set $S$ called the {\it musical space}, an {\it interval group} $\textit{IVLS}$, and an {\it interval function} $\textit{int}\colon S \times S \to \textit{IVLS}$ such that:
\begin{enumerate}
\item
for each $s \in S$ and $i \in \textit{IVLS}$, there is a unique $t \in S$ such that $\text{int}(s,t)=i$, and
\item
for all $s,t,u \in S$, $\textit{int}(s,t)\textit{int}(t,u)=\textit{int}(s,u)$.
\end{enumerate}
Lewin proved \cite[pages 157-158]{LewinGMIT} that this intervallic notion is equivalent to a transformational reformulation: a {\it simply transitive group action} of a group $\textit{SIMP}$ on a set $S$ is a group action in which each pair $s,t \in S$ admits a unique $g \in \textit{SIMP}$ such that $gs=t$. To obtain the transformational description from the intervallic description, one takes $\textit{SIMP}$ to be the group of transpositions. For an example of both perspectives, recall from \cite{kolman} that any group $G$ has an associated {\it canonical generalized interval system} $\mathbf{G}=(G,G,\textit{int}_G)$ with $\textit{int}_G(g,h)=g^{-1}h$. The
{\it transposition} $T_i \colon G \to G$ is then $T_i(g)=ig$, and $\textit{SIMP} \cong G^{\text{op}}$.

See Lewin \cite[Chapter 7]{LewinGMIT} for an explanation of conceptual differences between the
intervallic and transformational descriptions. See also Clampitt \cite{Clampitt} for an explanation of how
different transformationally described generalized interval systems suggest different hearing strategies. We have chosen to develop morphisms in both the intervallic and transformation descriptions in this article to have a self-contained, streamlined treatment of these complementary points of view, all in one place.

The equivalence of simply transitive group actions and affine spaces has been considered in a variety of mathematical contexts. See for instance \cite{Berger}.

A familiar example of a generalized interval system in the transformational formulation has set $S$ the set of major and minor chords, and transpositions the mod 12 $T/I$-group. Note that both transpositions {\it and inversions} are called transpositions in the GIS description of this example.

We develop morphisms in the intervallic picture first, starting with Definition~\ref{defn:intervallic_morphism}, and turn to the transformational picture later, starting with Definition~\ref{defn:transformational_morphism}.

\begin{defn}[Morphism of GIS's in intervallic description, Definition~9 of \cite{kolman}] \label{defn:intervallic_morphism}
A {\it morphism of generalized interval systems}
$$\xymatrix{(f,\varphi)\colon (S_1,\textit{IVLS}_1,\textit{int}_1) \ar[r] & (S_2,\textit{IVLS}_2,\textit{int}_2)}$$ consists of a function $f\colon S_1 \to S_2$ and a group homomorphism $\varphi\colon \textit{IVLS}_1 \to \textit{IVLS}_2$ such that
\begin{equation} \label{equ:intervallic_morphism}
\textit{int}_2\left(f(s),f(t)\right)=\varphi(\textit{int}_1(s,t))
\end{equation}
for all $s,t \in S_1$. That is, the diagram
$$\xymatrix@C=3pc{S_1 \times S_1 \ar[r]^-{f \times f} \ar[d]_{\textit{int}_1} & S_2 \times S_2 \ar[d]^{\textit{int}_2}
\\  \textit{IVLS}_1 \ar[r]_{\varphi} & \textit{IVLS}_2 }$$
commutes.
\end{defn}

Note that $\varphi$ does not determine $f$. For example, equation~\eqref{equ:intervallic_morphism} remains true if we replace $f$ by $p \circ f$ where $p\colon S_2 \to S_2$ is any interval preserving function on $S_2$.

We denote the category of generalized interval systems and their morphisms by $\mathbf{GIS}$. Composition and identity in $\mathbf{GIS}$ are $(g, \psi) \circ (f,\varphi)=(g \circ f, \psi \circ \varphi)$ and $1_{ (S,\textit{IVLS},\textit{int})}=(1_S,1_\textit{IVLS})$.

\begin{examp}[Interval preserving maps give rise to morphisms] \label{examp:interval_preserving_maps_give_rise_to_morphisms}
If $f \colon S \to S$ is an interval preserving function for $(S,\textit{IVLS},\textit{int})$, then $(f,\text{Id}_{\textit{IVLS}})$ is a morphism  $(S,\textit{IVLS},\textit{int}) \to (S,\textit{IVLS},\textit{int})$ with $\varphi=\text{Id}_{\textit{IVLS}}$. In fact, $(f,\text{Id}_{\textit{IVLS}})$ is automatically an {\it automorphism} by Proposition~\ref{prop:monic/epic_morphisms_for_intervallic}.
\end{examp}

\begin{examp} \label{examp:GISmorphism_from_group_homomorphism}
If $\varphi \colon G_1 \to G_2$ is any group homomorphism and $a \in G_2$, then $(f,\varphi):=(a\varphi,\varphi)$ is a morphism of the associated canonical generalized interval systems $\mathbf{G_1} \to \mathbf{G_2}$, since
$$\textit{int}_2(f(s),f(t))=f(s)^{-1}f(t)=\varphi(s)^{-1}(a^{-1}a)\varphi(t)=\varphi(s^{-1}t)=\varphi(\textit{int}_1(s,t)).$$
\end{examp}

In any category, it is useful to be able to recognize embeddings, quotient maps, and isomorphisms. One can show that embeddings and quotients in $\mathbf{GIS}$ are monic and epic morphisms, and that a morphism is an isomorphism if and only if it is monic and epic. Proposition~\ref{prop:monic/epic_morphisms_for_intervallic} states that
to determine whether a morphism of generalized interval systems is monic/epic/iso, it suffices to consider whether the musical space map $f$ or the interval map $\varphi$ is injective/surjective/bijective. We make use of the upcoming Proposition~\ref{prop:monic/epic_morphisms_for_intervallic} several times in this paper, notably in  Example~\ref{examp:interval_preserving_maps_give_rise_to_morphisms},
Proposition~\ref{prop:GrpAff_equivalent_to_GIS}, Example ~\ref{examp:centralizer_elements_are_morphisms}, Proposition~\ref{prop:transfer_of_group_actions},  and at the end of Section~\ref{sec:Octatonic_Example}.

Recall that a morphism $m\colon b \to c$ in a category $\mathcal{C}$ is said to be {\it monic} if for any two morphisms $f_1, f_2 \colon a \to b$ in $\mathcal{C}$, we conclude from $m\circ f_1= m \circ f_2$ that $f_1=f_2$. A morphism $p\colon b \to c$ in $\mathcal{C}$ is said to be {\it epic} if for any two morphisms $g_1, g_2 \colon b \to c$ in $\mathcal{C}$, we conclude from $g_1 \circ p  = g_2 \circ p$ that $g_1=g_2$, see \cite[page 19]{MacLaneWorking}.

\begin{prop}[Monic/Epic/Iso morphisms, intervallic picture] \label{prop:monic/epic_morphisms_for_intervallic}
Let $(f,\varphi)$ be a morphism of generalized interval systems as in Definition~\ref{defn:intervallic_morphism}.
\begin{enumerate}
\item
The following are equivalent.
\begin{enumerate} \label{item:monic}
\item
The morphism $(f,\varphi)$ is monic.
\item
The function $f$ is injective. \label{item:item:monic_f}
\item
The homomorphism $\varphi$ is injective. \label{item:item:monic_varphi}
\end{enumerate}
\item \label{item:epic}
The following are equivalent.
\begin{enumerate}
\item
The morphism $(f,\varphi)$ is epic.
\item
The function $f$ is surjective. \label{item:item:epic_f}
\item
The homomorphism $\varphi$ is surjective. \label{item:item:epic_varphi}
\end{enumerate}
\item \label{item:iso}
The following are equivalent.
\begin{enumerate}
\item
The morphism $(f,\varphi)$ is an isomorphism.
\item
The function $f$ is bijective. \label{item:item:iso_f}
\item
The homomorphism $\varphi$ is bijective. \label{item:item:iso_varphi}
\end{enumerate}
\end{enumerate}
\end{prop}
\begin{pf}
Since the composition $(g,\psi) \circ (f,\varphi)$ is simply $(g\circ f,\psi \circ \varphi)$ we immediately see that
$(f,\varphi)$ is monic/epic/iso if and only if both $f$ and $\varphi$ are monic/epic/iso. But a function or group homomorphism is monic/epic/iso if and only if it is injective/surjective/bijective.

For \ref{item:item:monic_f}$\Rightarrow$\ref{item:item:monic_varphi}, suppose $f$ is injective and $\varphi(i)=\varphi(j)$. Then there exist $s,t,u$ in $S_1$ such that $\textit{int}_1(s,t)=i$ and $\textit{int}_1(s,u)=j$. By equation \eqref{equ:intervallic_morphism} we have
\[
\textit{int}_2(f(s),f(t))=\varphi(\textit{int}_1(s,t))=\varphi(\textit{int}_1(s,u))=\textit{int}_2(f(s),f(u))
\]
so that $f(t)=f(u)$, then $t=u$, and finally $i=j$.

For \ref{item:item:monic_varphi}$\Rightarrow$\ref{item:item:monic_f}, suppose $\varphi$ is injective and $f(s)=f(t)$. Then equation \eqref{equ:intervallic_morphism} implies $\textit{int}_1(s,t)=e_{\textit{IVLS}_1}$ by the injectivity of $\varphi$, so that $s=t$. This completes the proof of claim \ref{item:monic}.

For \ref{item:item:epic_f}$\Rightarrow$\ref{item:item:epic_varphi}, suppose $f$ is surjective and $k \in \textit{IVLS}_2$. Then there exist $s',t' \in S_2$ such that $k=\textit{int}_2(s',t')$, and by equation \eqref{equ:intervallic_morphism} and the surjectivity of $f$ we have $s,t \in S_1$ with
\[
k=\textit{int}_2(s',t')=\textit{int}_2(f(s),f(t))=\varphi(\textit{int}_2(s,t)).
\]

For \ref{item:item:epic_varphi}$\Rightarrow$\ref{item:item:epic_f}, suppose $\varphi$ is surjective and $t' \in S_2$. Fix $s \in S_1$. Then for some $i \in \textit{IVLS}_1$ and some $t \in S_1$ we have
\[
\textit{int}_2(f(s),t')=\varphi(i)=\varphi(\textit{int}_1(s,t))=\textit{int}_2(f(s),f(t)),
\]
which implies $t'=f(t)$. This completes the proof of claim \ref{item:epic}.

The equivalences in \ref{item:iso} follow from \ref{item:monic} and \ref{item:epic} and the fact that a group homomorphism (or function) is an isomorphism if and only if it is bijective.
\end{pf}

As an example of Proposition~\ref{prop:monic/epic_morphisms_for_intervallic}, if $\varphi \colon G_1 \to G_2$ is a group homomorphism, then its associated morphism $(\varphi,\varphi)$ of canonical generalized intervals systems $\mathbf{G_1} \to \mathbf{G_2}$ in Example~\ref{examp:GISmorphism_from_group_homomorphism} is monic if and only if $\varphi$ is injective, and epic if and only if $\varphi$ is surjective. Thus, morphisms of generalized interval systems are in general neither monic nor epic. For instance, the quotient map $\mathbb{Z} \to \mathbb{Z}_n$ induces a surjective map between the associated canonical generalized interval systems. If $G_0$ is a subgroup of $G$, then the inclusion $G_0 \hookrightarrow G$ induces an inclusion of the associated canonical generalized interval systems.

We next turn to morphisms in the transformational picture of generalized interval systems, namely morphisms of simply transitive group actions.
For the transformational picture, we first introduce the more general notion of a morphism of group actions, and then turn to the more specific actions of simply transitive group actions, since the actions in some applications are not simply transitive.

\begin{defn}[Morphisms of group actions and simply transitive group actions] \label{defn:transformational_morphism}
Suppose $(G_1,S_1)$ and $(G_2,S_2)$ are group actions, not necessarily simply transitive.
A {\it morphism of group actions}\footnote{In order to have the same notation as GIS morphisms, we have chosen to put the function $f\colon S_1\to S_2$ before the group homomorphism $\varphi\co G_1 \to G_2$ in our notation $(f,\varphi)$ for a morphism of group actions, even though the notation for group actions $(G_1,S_1)$ puts the group before the set.}
$$\xymatrix{(f,\varphi)\colon (G_1,S_1) \ar[r] & (G_2,S_2)}$$ consists of a function $f\colon S_1 \to S_2$ and a group homomorphism $\varphi \colon G_1 \to G_2$ such that
$$f(gs)=\varphi(g)f(s)$$
for all $g \in G_1$ and all $s \in S_1$. That is, the diagram
$$\xymatrix@C=3pc{G_1 \times S_1 \ar[r]^{\varphi \times f} \ar[d]_{\text{action}} & G_2 \times S_2 \ar[d]^{\text{action}} \\ S_1 \ar[r]_f & S_2 }$$
commutes. If the actions $(G_1,S_1)$ and $(G_2,S_2)$ are simply transitive, then the group action morphism $(f,\varphi)$ is said to be a {\it morphism of simply transitive group actions}.
\end{defn}

We denote the category of simply transitive group actions and their morphisms by $\mathbf{SimpTransGrpAct}$. Composition and identity in this category are $(g, \psi) \circ (f,\varphi)=(g \circ f, \psi \circ \varphi)$ and $1_{ (G,S)}=(1_S,1_G)$.

The transformational analogue of Proposition~\ref{prop:monic/epic_morphisms_for_intervallic} for monic/epic/iso morphisms of simply transitive group actions clearly holds.

Our main examples of morphisms of simply transitive group actions arise from affine maps in the situation of groups generated by contextual inversions, see Corollary~\ref{cor:affines_commute_with_subgroups_generated_by_contextual_inversions}.
In the Schoenberg analysis in Section~\ref{sec:Analysis_of_Schoenberg}, the affine map  ${}^77$ provides us with a morphism to model a piecewide narrative. We present two more examples in the following Example and Proposition.

\begin{examp}[Commuting elements give rise to morphisms] \label{examp:centralizer_elements_are_morphisms}
If $G\leq \text{Sym}(S)$ acts simply transitively on $S$, and if $f$ is in the centralizer $C_{\text{\rm Sym}(S)}(G)$, then $(f,\text{Id}_G)$ is a morphism $(G,S) \to (G,S)$. In fact, $(f,\text{Id}_G)$ is even an automorphism by the transformational analogue of Proposition~\ref{prop:monic/epic_morphisms_for_intervallic}. This example is the transformational analogue of Example~\ref{examp:interval_preserving_maps_give_rise_to_morphisms}.
\end{examp}

\begin{prop}[Transfer of group actions] \label{prop:transfer_of_group_actions}
Suppose $G\leq \text{\rm Sym}(S_1)$ acts  on $S_1$, and $f\colon S_1 \to S_2$ is a bijection. Then the group $fGf^{-1} \leq \text{\rm Sym}(S_2)$ acts on $S_2$ and $(f,\gamma^f)\colon (G_1,S_1) \to (fGf^{-1},S_2)$ is an isomorphism of group actions, where $\gamma^f(g)=fgf^{-1}$. If $G$ acts simply transitively, then so does $fGf^{-1}$, and its dual is the $f$-conjugate of the dual to $G$.
\end{prop}
\begin{pf}
To verify $(f, \gamma_f)$ is a morphism, the diagram in Definition~\ref{defn:transformational_morphism} can be checked directly. Since $f$ is a bijection, the morphism $(f, \gamma_f)$ is an isomorphism.
Simple transitivity of $fGf^{-1}$ can be directly checked, while the claim about duals follows from the fact $\gamma_f\co \mathrm{Sym}(S_1) \to \mathrm{Sym}(S_2)$ is an isomorphism.
\end{pf}

\begin{examp}
As an example of Proposition~\ref{prop:transfer_of_group_actions}, if $f\co \mathbb{Z}_{12} \to \mathbb{Z}_{12}$ is a bijection, and $S$ is
a set of $n$-tuples in $\mathbb{Z}_{12}$ upon which a group $G$ acts simply transitively, then, by computing componentwise, $f$ induces a bijection $S \to f(S)$ (also called $f$) and $f G f^{-1}$ acts simply transitively on $f(S)$, with its dual being the conjugate of the dual of $G$. In particular, this is the relationship between the first and second rows of Figure~\ref{fig:dual_groups_to_G0}  with $f$ equal to multiplication by $7$ (which commutes with the $H$-groups).  The rows of Figure~\ref{fig:dual_groups_to_G0} are also related to the respective rows of Figure~\ref{fig:dual_groups_to_T1G0T1Inverse} in this way, where $f=T_1$.
\end{examp}

The rest of this section is dedicated to proving the equivalence of three categories: $\mathbf{GIS}$ (generalized interval systems and their morphisms),
$\mathbf{SimpTransGrpAct}$ (simply transitive group actions and their morphisms), and $\mathbf{GrpAff}$ (groups and affine morphisms). It has been known since at least Lewin's monograph \cite[Chapter 7]{LewinGMIT} that generalized interval systems and simply transitive group actions are ``equivalent''. For the sake of completeness we quickly derive an equivalence of categories $\mathbf{GIS}\to\mathbf{SimpTransGrpAct}$, and an equivalence of categories $\mathbf{GrpAff} \to \mathbf{GIS}$.

\begin{prop}
The category of generalized interval systems is equivalent to the category of simply transitive group actions via the functor $$\xymatrix{F\colon \mathbf{GIS} \ar[r] & \mathbf{SimpTransGrpAct}}$$
$$\xymatrix{(S,\textit{IVLS},\textit{int}) \ar@{|->}[r] & (\textit{SIMP},S)}$$
$$\xymatrix{(f,\varphi) \ar@{|->}[r] & (f,\overline{\varphi}),}$$
where $\overline{\varphi}(T_i)=T_{\varphi(i)}$.
\end{prop}
\begin{pf}
For the details of the functor on the object level, see Lewin \cite[pages 157--158]{LewinGMIT}.

If $(f,\varphi)$ is a GIS morphism, then $F(f,\varphi)=(f,\overline{\varphi})$ is a morphism of simply transitive group actions because $\textit{int}_1(s,T_i(s))=i$ implies $\textit{int}_2(f(s),f\left(T_i(s)\right))=\varphi(i)$ and therefore $f\left(T_i(s)\right)=T_{\varphi(i)}\left(f(s)\right),$
so that $f\left(T_i(s)\right)=\overline{\varphi}(T_i)f(s)$ as required in Definition~\ref{defn:transformational_morphism}. The function $\overline{\varphi}$ is a group homomorphism because
\begin{equation} \label{equ:bar_is_homomorphism}
\overline{\varphi}(T_iT_j)=\overline{\varphi}(T_{ji})=T_{\varphi(ji)}=T_{\varphi(j)\varphi(i)}=T_{\varphi(i)}T_{\varphi(j)}=\overline{\varphi}(T_i)\overline{\varphi}(T_j).
\end{equation}

To prove that $F$ is an equivalence of categories, it suffices to prove $F$ is faithful, full, and essentially surjective \cite[Theorem 1, page 93]{MacLaneWorking}.

For faithfulness, if $(f,\overline{\varphi})=(g,\overline{\psi})$, then $f=g$ and $T_{\varphi(i)}=T_{\psi(i)}$ for all $i \in \textit{IVLS}_1$, which implies $\varphi(i)=\psi(i)$ for all $i \in \textit{IVLS}_1$ because $j \mapsto T_j$ is an anti-isomorphism $\textit{IVLS}_2\to \textit{SIMP}_2$.

For fullness, if $(f,\nu)\co (\textit{SIMP}_1,S_1) \to (\textit{SIMP}_2,S_2)$, then define $$\varphi\co \textit{IVLS}_1 \to \textit{IVLS}_2$$ via the equation $\nu(T_i)=T_{\varphi(i)}$ (this uniquely defines $\varphi$ because of the anti-isomorphism mentioned above for $\textit{SIMP}_1$ and $\textit{SIMP}_2$). Then $\overline{\varphi}=\nu$ and $\varphi$ is a group homomorphism by reversing the argement of equation \eqref{equ:bar_is_homomorphism}.

For essential surjectivity, if $(G,S)$ is a simply transitive group action, define a generalized interval system with musical space $S$, interval group $G^{\text{op}}$, and $\textit{int}(s,t):=$ that $g\in G$ such that $gs=t$. Then $F(S,\textit{IVLS},\textit{int})$ is isomorphic to $(G,S)$ using the isomorphism $\textit{SIMP}\cong (G^\text{op})^\text{op}\cong G$ and the identity on $S$.
\end{pf}

An equivalence between groups with affine maps and generalized interval systems is essentially contained in Kolman's paper \cite{kolman}; we give a proof next.

Let $\mathbf{GrpAff}$ denote the category with groups as objects and with affine maps of groups as morphisms. A function $f\co G \to H$ between groups is an {\it affine map of groups} if it is a group homomorphism followed by a left multiplication, that is, if it has the form $f(g)=a\varphi(g)$ for some $a \in H$ and some group homomorphism $\varphi\co G \to H$. In this case, $a=f(e)$, so $a$ and $\varphi$ are uniquely determined by $f$.

\begin{prop} \label{prop:GrpAff_equivalent_to_GIS}
The category of groups with affine maps is equivalent to the category of generalized interval systems via the ``canonical generalized interval system'' functor
$$\xymatrix{C \colon \mathbf{GrpAff} \ar[r] & \mathbf{GIS}}$$
$$\xymatrix{G \ar@{|->}[r] & \mathbf{G}=(G,G,\textit{int}_G)}$$
$$\xymatrix{a\varphi \ar@{|->}[r] & (a\varphi,\varphi)}.$$
\end{prop}
\begin{pf}
If $C(a_1\varphi_1,\varphi_1)=C(a_2\varphi_2,\varphi_2)$, then $\varphi_1=\varphi_2$, and evaluation at $e_G$ gives $a_1=a_2$, so that $C$ is faithful. The functor $C$ is also full, for if $(f,\varphi)\colon \mathbf{G} \to \mathbf{H}$ is a morphism of (canonical) generalized interval systems, then taking $s=e$ in equation~\eqref{equ:intervallic_morphism} yields $f(t)=f(e)\varphi(t)$, so $f(e)\varphi$ is a preimage of $(f,\varphi)$.

For the essential surjectivity, consider any GIS $(S,\textit{IVLS},\textit{int})$. Let $G:=\textit{IVLS}$ and fix an element $s_0 \in S$.  We define an isomorphism $(f,\varphi)\colon \mathbf{G} \to  (S,\textit{IVLS},\textit{int})$ by
$f(g):=T_gs_0$ and $\varphi:=\text{id}_G$. The pair $(f,\varphi)$ is a morphism because
\[
\textit{int}(f(g),f(h))=\textit{int}(T_gs_0,T_hs_0)=\textit{int}(T_gs_0,s_0)\textit{int}(s_0,T_hs_0)=g^{-1}h.
\]
The pair $(f,\varphi)$ is an isomorphism because $\varphi$ is a bijection (see Proposition~\ref{prop:monic/epic_morphisms_for_intervallic}).
\end{pf}

Proposition~\ref{prop:GrpAff_equivalent_to_GIS} illustrates one reason why Mazzola and his collaborators \cite{MazzolaToposOfMusic, MazzolaGeometrieDerToene, MazzolaGruppenUndKategorien} use affine maps rather than mere group homomorphisms.

\section{Sub Dual Groups and Contextual Groups} \label{sec:Sub_Dual_Groups_and_Contextual_Groups}

Let $\text{\rm Sym}(S)$ be the symmetric group on the set $S$. Two subgroups $G$ and $H$ of the symmetric group $\text{\rm Sym}(S)$ are {\it dual in the sense of Lewin} \cite[page 253]{LewinGMIT} if their natural actions on $S$ are simply transitive and each is the centralizer of the other, that is,
$$C_{\text{\rm Sym}(S)}(G)=H \;\; \text{ and } \;\; C_{\text{\rm Sym}(S)}(H)=G.$$ Each generalized interval system gives rise to a dual pair: its group of transpositions and its group of interval preserving bijections, see \cite[pages 251--253]{LewinGMIT}. If subgroups $G_0$ and $H_0$ of dual groups $G$ and $H$ preserve a subset $S_0 \subseteq S$, and $G_0$ and $H_0$ are dual in $\text{\rm Sym}(S_0)$, then $G_0$ and $H_0$ are said to be {\it sub dual groups}\footnote{Motivation for the notion of sub dual groups can be found in Figures 2.12 and 2.13 on page 34 of John Clough's article \cite{cloughFlipFlop}. There he considers sub dual groups by superimposing their Cayley graphs, though he does not explicitly thematize the concept of sub dual groups.} of $G$ and $H$.

A simple method for constructing sub dual groups from dual groups was given by Fiore--Noll in \cite{fiorenoll2011}.

\begin{thm}[Construction of sub dual groups, Theorem 3.1 of \cite{fiorenoll2011}] \label{thm:sub_dual_groups}
Let $G,H \leq \text{\rm Sym}(S)$ be dual groups, $G_0$ a subgroup of $G$, and $s_0$ an element of $S$. Let $S_0$ be the orbit of $s_0$ under the action of $G_0$. Then the following hold.
\begin{enumerate}
\item \label{thm:sub_dual_groups:G_0_simple_transitivity}
The group $G_0$ acts simply transitively on $S_0$.
\item \label{thm:sub_dual_groups:s_0_determines_G_0}
If $g \in G$ and $gs_0$ is in $S_0$, then $g \in G_0$. In particular, if $g \in G$ and $gs_0$ is in $S_0$, then $g$ preserves $S_0$ as a set, that is, $g(S_0) \subseteq S_0$.
\item \label{thm:sub_dual_groups:H_0_simple_transitivity}
Let $H_0$ denote the subgroup of $H$ consisting of those elements $h\in H$ with $hs_0 \in S_0$.  Then $H_0$ acts simply transitively on $S_0$.
\item
Restriction from $S$ to $S_0$ embeds $G_0$ and $H_0$ in $\text{\rm Sym}(S_0)$. We denote their images in $\text{\rm Sym}(S_0)$ by $G_0\vert_{S_0}$ and $H_0\vert_{S_0}$.
\item \label{thm:sub_dual_groups:extensions}
If $g \in \text{\rm Sym}(S_0)$ and $g$ commutes with $H_0\vert_{S_0}$, then $g$ admits a unique extension to $S$ which belongs to $G$. This extension necessarily belongs to $G_0$ by \ref{thm:sub_dual_groups:s_0_determines_G_0}. Similarly, if $h \in \text{\rm Sym}(S_0)$ commutes with $G_0\vert_{S_0}$, then $h$ admits a unique extension to $S$ which belongs to $H$. This extension necessarily belongs to $H_0$.
\item \label{thm:sub_dual_groups:restricted_subgroups_dual}
The groups $G_0\vert_{S_0}$ and $H_0\vert_{S_0}$ are dual in $\text{\rm Sym}(S_0)$.
\end{enumerate}
\end{thm}
If we consider $G_0$ acting on another element $ks_0$ in $S$, with $k \in H$, then the orbit of $ks_0$ is $kS_0$ and the associated dual group is $kH_0k^{-1}$ \cite[Corollary 3.3]{fiorenoll2011}. The inclusions $(G_0,S_0) \hookrightarrow (G,S)$ and $(H_0,S_0) \hookrightarrow (H,S)$ are monic morphisms of simply transitive group actions.

Dual groups $G$ and $H$ useful for constructing sub dual groups, especially for our octatonic example, are given by the $T/I$-group and the generalized contextual group of Fiore--Satyendra \cite{fioresatyendra2005}, which we now recall. For another treatment of transformations in serial contexts, see Hook \cite{hookUTT2002} and Hook--Douthett \cite{hook2008uniform} in terms of uniform triadic transformations.\footnote{The article \cite{fioresatyendra2005} of Fiore--Satyendra contains some comparisons with and comments on the uniform triadic transformations of Hook \cite{hookUTT2002}. Particular points of intersection of the present approach with Hook--Douthett \cite{hook2008uniform},  such as the RICH transformation, are discussed in \cite{fiorenollsatyendraMCM2013}.
}

Suppose $X = (x_1,\dots, x_n)$ is a pitch-class segment\footnote{A {\it pitch-class segment} is an ordered subset of $\mathbb{Z}_m$. We use parentheses to denote an ordered subset of $\mathbb{Z}_m$ as an $n$-tuple: $( x_1, \dots , x_n )$. We do not use the traditional musical notation $\langle x_1, \dots , x_n \rangle$ for pitch-class segments because it clashes with the mathematical notation for the subgroup generated by $x_1, \dots, x_n$, which we will also need on occasion.

The contextual groups of Fiore--Satyendra \cite{fioresatyendra2005} use {\it ordered} pitch-class sets (=pitch-class segments) as opposed to {\it unordered} pitch-class sets for a variety reasons. Musically speaking, the ordering allows an extension of contextual inversions to arbitrary pitch-class segments, even if the underlying pitch-class set has a symmetry.  Mathematically speaking, the ordering enables a simple definition of contextual inversion: for instance, $K(Y):= I_{y_1 + y_2}(Y)$ is easily seen to be a well-defined operation which sends $Y$ to the inversion of $Y$ which has the entries $y_1$ and $y_2$ switched. Without an ordering, the definition would be clumsy, perhaps referring to some artificial root. The ordering of the initial pitch-class segment $X$ does not presume a temporal ordering of the notes in the work of music at hand. For instance, if $X=(0,4,7)$ then there is clearly a bijection from the $T/I$-class of $(0,4,7)$ to the $T/I$-class of $\{0,4,7\}$ given by $(0,4,7) \mapsto \{0,4,7\}$ et cetera. In the ordered setting of $(0,4,7)$, we define $P$ as $P(Y)=I_{y_1+y_3}(Y)$ and we see that this coincides (via the aforementioned bijection) with the usual definition of $P$. See  \cite[Footnote 20]{fioresatyendra2005}.

However, in \cite{fiorenollsatyendraMCM2013} we use permutations to show how to extend $P$, $L$, and $R$ to the various reordings of major and minor chords, so that orderings of pitch-class segments {\it do indeed coincide} with the temporal ordering of the notes in a given work of music. Notably, $P$, $L$, and $R$ are {\it not} globally defined on reorderings by contextual inversion formulas such as $P(Y)=I_{y_1+y_3}(Y)$, rather they are defined via conjugations. See for instance the incorporation of the retrograde at the end of Section~\ref{sec:Analysis_of_Schoenberg}. } with entries in $\mathbb{Z}_{m}$. Let $S$ be the orbit of $X$ under the componentwise action of the mod $m$ $T/I$-group, in other words, $S$ is the serialized set class of $X$. The orbit $S$ has two types of elements: reflections of $X$ and inversions of $X$. A pitch-class segment $Y=( y_1,\dots, y_n ) \in S$ is said to be a {\it $T$-form} if it is a translate of $X$, that is, if there exists some $i \in \mathbb{Z}_m$ such that $Y=T_i(X)$. Similarly, $Y$ is said to be an {\it $I$-form} if it is an inversion of $X$, that is, if there exists some $j \in \mathbb{Z}_m$ such that $Y=I_j(X)$. Let $K,Q_i\colon S \to S$ be the bijections defined on a pitch-class segment $Y=( y_1,\dots, y_n ) \in S$ as follows.
$$
\begin{array}{rl}
K(Y) &:= I_{y_1 + y_2}(Y) \\
& \\
Q_i(Y) &:= \left\{
\begin{array}{ll}
T_i(Y) & \text{if $Y$ is a $T$-form of $X$} \\
T_{-i}(Y) & \text{if $Y$ is an $I$-form of $X$.}
\end{array} \right.
\end{array}
$$

The subgroup of $\text{Sym}(S)$ generated by $K$ and $Q_1$ is called in \cite{fioresatyendra2005} the {\it generalized contextual group associated to $X$}, or simply the {\it contextual group}. Suppose now that $X$ satisfies the {\it tritone condition}\footnote{The {\it tritone condition} requires that there exist two distinct pitch classes $x_q,x_r$ in $X$ which span an interval other than $m/2$.}. Then the $T/I$-group and the contextual group act simply transitively on $S$, are dual, and are both dihedral of order $2m$, as proved in \cite[Theorem 4.2 and Corollary 4.3]{fioresatyendra2005}. Each element is either of the form $Q_i$ or $Q_iK$. If $J^{q,r}$ is defined to exchange $y_q$ and $y_r$, rather than $y_1$ and $y_2$ as above in $K$, then the group generated by $J^{q,r}$ and $Q_1$ is equal to the generalized contextual group \cite[Corollary 4.4]{fioresatyendra2005}. In particular, each bijection $J^{q,r}$ defined in terms of common-tone preservation
\begin{equation} \label{equ:contextual_inversion}
J^{q,r}\colon Y \mapsto I_{y_q+y_r}Y
\end{equation}
is an element of the contextual group acting on the $T/I$-orbit of $X$. Such bijections are called {\it contextual inversions}.\footnote{Related work is Childs \cite{childs}, Gollin \cite{gollin}, and Kochavi \cite{kochavi}.} We remark here that $J^{q,r}$ is not in the dual group to the $T/I$-group with permutations $\Sigma_n$ when $S$ is enlarged to include all possible orderings of major and minor triads.

If $X=( 0, 4, 7 )$ in $\mathbb{Z}_{12}$, then the associated contextual group is precisely the neo-Riemannian $PLR$-group, and we have on the $T/I$-orbit of $(0,4,7)$
the following equalities.\footnote{Beyond the root position of major and minor triads, i.e. on reorderings of the $T/I$-orbit of $(0,4,7)$, we no longer have $P$, $L$, and $R$ equal to $J^{1,3}$, $J^{2,3}$, and $J^{1,2}$. Rather, $P$, $L$, and $R$ are extended to other orderings of major and minor triads via conjugation with the relevant permutations. See for instance the incorporation of the retrograde at the end of Section~\ref{sec:Analysis_of_Schoenberg}. The relationship between $\mathrm{RICH}$-type transformations, contextual inversions, the $PLR$-group, and voice permutations is clarified in \cite{fiorenollsatyendraMCM2013}.}
$$\aligned
P(y_1,y_2,y_3)&=J^{1,3}(Y) &=I_{y_1+y_3}(y_1,y_2,y_3)&=(y_3,-y_2+y_1+y_3,y_1) \\
L(y_1,y_2,y_3)&=J^{2,3}(Y) &=I_{y_2+y_3}(y_1,y_2,y_3)&=(-y_1+y_2+y_3,y_3,y_2) \\
R(y_1,y_2,y_3)&=J^{1,2}(Y) &=I_{y_1+y_2}(y_1,y_2,y_3)&=(y_2,y_1,-y_3+y_1+y_2)
\endaligned$$
The neo-Riemannian $PLR$-group is the same as the {\it Schritt-Wechsel} group. The name  ``{\it Schritt-Wechsel} group'' is often used for the neo-Riemannian $PLR$-group in the literature because each operation $Q_i$ is called a {\it Schritt} while each operation $Q_iK$ is called a {\it Wechsel}.

If $X$ is a non-consonant triad satisfying the tritone condition, then the $P$, $L$, $R$ analogues $J^{1,3}$, $J^{2,3}$, and $J^{1,2}$ on $S$
are in the contextual group associated to $X$, though they may not generate it.\footnote{For example, consider the pitch-class segment $X=(0,4,10)$, as we will in Section~\ref{sec:Octatonic_Example}. Since the entries of $X$ are all of the same parity (even), the $P$, $L$, and $R$ operations in equation \eqref{equ:contextual_inversion} will also preserve parity. Therefore, no composite of $P$, $L$, and $R$ maps $(0,4,10)$ to $(1,5,11)$. The subgroup $\langle P,L,R \rangle$ of the contextual group does not act simply transitively, so must be properly contained in the (simply transitive) contextual group. Catanzaro first observed geometrically that the generalized $PLR$-group for a given pitch-class set may not act simply transitively in \cite{Catanzaro}: the {\it Tonnetz} has multiple components when the generating pitch-class set is entirely even. Though Catanzaro works with unordered pitch-class sets (as opposed to the ordered pitch-class segments in the contextual group of \cite{fioresatyendra2005}), the example $(0,4,10)$ is an example for both, as it has no symmetries.}

Contextual inversions $J^{q,r}$, as defined in equation \eqref{equ:contextual_inversion}, are compatible with affine maps, as we will need in Section~\ref{sec:Analysis_of_Schoenberg} for our analysis of Schoenberg's String Quartet in $D$ minor.

\begin{thm}[Contextual inversions commute with all affine maps] \label{thm:contextual_inversions_commute_with_affine_maps}
If $f\colon \mathbb{Z}_m \to \mathbb{Z}_m$ is an affine map, and $J^{q,r}$ is the contextual inversion in \eqref{equ:contextual_inversion}, then $f\circ J^{q,r}(Y)=J^{q,r} \circ f(Y)$ for any pitch-class segment $Y$.
\end{thm}
\begin{pf}
Suppose $f\colon \mathbb{Z}_m \to \mathbb{Z}_m$ is $f(z)=az+b$ for fixed $a,b \in \mathbb{Z}_m$ and
let $y_i$ be the $i$-th component of $Y$. Then
$$
\aligned
\left( f\circ J^{q,r}(Y) \right)_i &= aI_{y_q+y_r}(y_i)+b \\
&= a(-y_i+y_q+y_r)+b \\
&= (-ay_i-b) +(ay_q+b) + (ay_r+b) \\
&= I_{(ay_q+b)+(ay_r+b)}(ay_i+b) \\
&= J^{q,r}\circ f(y_i).
\endaligned
$$
\end{pf}

We now have a morphism of (simply transitive) group actions  as described in Section~\ref{sec:Morphisms}.

\begin{cor}[Affine maps give rise to morphisms of group actions when the domain group is generated by contextual inversions] \label{cor:affines_commute_with_subgroups_generated_by_contextual_inversions}
Let $X$ be a pitch-class segment in $\mathbb{Z}_m$, $S$ its $T/I$-orbit, $f\colon \mathbb{Z}_m \to \mathbb{Z}_m$ an affine map, and $G_0$ a subgroup of the contextual group for $X$ generated by contextual inversions. Let $G'$ and $S'$ be the contextual group and $T/I$-orbit for $f(X)$, and $\varphi\colon G_0 \to G'$ the homomorphism $J^{q,r} \mapsto J^{q,r}$.  Then
$$\xymatrix{(f,\varphi) \colon (G_0, S) \ar[r] & (G',S')}$$
is a morphism of group actions. The restriction to $S_0:=G_0X$ is a morphism
$$\xymatrix{(f,\varphi) \colon (G_0, S_0) \ar[r] & (G',S')}$$
of simply transitive group actions as in Definition~\ref{defn:transformational_morphism}.
\end{cor}
\begin{pf}
The map $f$ commutes with the generators of $G_0$ by Theorem~\ref{thm:contextual_inversions_commute_with_affine_maps}, so $f$ commutes with all of $G_0$.
\end{pf}

\begin{rmk} \label{rmk:suggestion_to_define_Qoverline} Corollary~\ref{cor:affines_commute_with_subgroups_generated_by_contextual_inversions} does not imply that affine maps commute with the $Q_i$ operations because the relations between $Q_i$ and contextual inversions depend on the respective pitch-class segment. For example, consider $X=(0,4,7)$, $f(z)=7z$, and $f(X)=(0,4,1)$. Then $L^{(0,4,7)}R^{(0,4,7)}=Q_5$ while $L^{(0,4,1)}R^{(0,4,1)}=Q_{11}$. So, commutation of $f$ with $L$ and $R$ implies $f Q_5 = Q_{11} f$.
\end{rmk}

Remark~\ref{rmk:suggestion_to_define_Qoverline} suggests that we define operations $\overline{Q}_i$ on the entire set of major, minor, jet, and shark chords in such a way that $\overline{Q}_i$ commutes with $M_7$ and $M_5$. Namely, we define
$$\overline{Q}_i(\text{major/minor chord})=Q_i(\text{major/minor chord})$$
$$\overline{Q}_i(\text{jet/shark chord})=Q_{7i}(\text{jet/shark chord})$$
so that $M_7\overline{Q}_i=\overline{Q}_iM_7$ and $M_5\overline{Q}_i=\overline{Q}_iM_5$ for all $i \in \mathbb{Z}_{12}$. In Remark~\ref{rmk:suggestion_to_define_Qoverline} we now have
$f\overline{Q}_5=\overline{Q}_5f$.

In fact, $\overline{Q}_i$ is part of the dual group to the group $\langle T/I, M_5, M_7 \rangle$ which acts simply transitively on its 48-element orbit of $(0,4,7)$, that is, on the entire set of major, minor, jet, and shark chords. More generally, we now describe the dual group to the full affine group
$$\mathrm{Aff}^\ast(\mathbb{Z}_{m})=\{x \mapsto ax+b\;\vert\; \text{$a$ is invertible in $\mathbb{Z}_m$}\}.$$

Suppose $X = (x_1,\dots, x_n)$ is a pitch-class segment with entries in $\mathbb{Z}_{m}$. Let $S$ be the orbit of $X$ under the componentwise action of the full affine group $\mathrm{Aff}^\ast(\mathbb{Z}_{m})$. For a pitch-class segment $Y=( y_1,\dots, y_n ) \in S$ and any $a \in \mathbb{Z}_m$ we have {\it Side Transformation (generalized Seitenwechsel)}
$$W_a((y_1, y_2, \dots, y_m)) := T_{(1-a) y_1} M_a ((y_1, y_2, \dots, y_m)).$$
Traditionally, the {\it Seitenwechsel} in Riemann (or {\it antinomic Wechsel} in Oettingen \cite[page 142]{Oettingen}) maps a major triad $(t, t + M3, t + P5)$ with tonic root $t$ into a minor triad $(t, t - M3, t - P5)$ with the same tone $t$ as its phonic root. In the Neo-Riemannian remake the transformation is the defined in the pitch class domain $\mathbb{Z}_{12}$ via $$\begin{array}{lll}W_{-1}((t, t+4, t+7)) & = & T_{2 t} \circ I ((t, t+4, t+7))\\ & = & T_{2 t} ((-t, -t-4, -t-7))\\ & = & (t, t-4, t-7)\end{array}$$

The definition of generalized Seitenwechsel $W_a$ assumes that the first coordinate of a pitch-class segment $Y=( y_1,\dots, y_n ) \in S$  is the generalized Oettingen-Root and should be preserved under the transformation. In other words, the Side Transformation is the unique transposition $T_k(M_a(Y))$ of the linear transform $M_a(Y)$ of $Y$, whose root coincides with the root of $Y$.

\noindent The Side Transformation commutes with any affine transformation  $f: \mathbb{Z}_{m} \to \mathbb{Z}_{m}$, say $f(y) = c y + d \,\, \mbox{mod} \,\, m$:
$$\begin{array}{lll}W_a(f(y_i)) & = & T_{(1-a) (c y_1 + d)} M_a (c y_i + d)\\
& = & a c y_i +  a d + (1-a) (c y_1 + d)\\
& = & c (a y_i + (1-a) y_1) + d \\
& = & f( W_a(y_i))
\end{array}$$

Furthermore we have the {\it generalized Schritt-Transformations} $\overline{Q}_j, \, j \in \mathbb{Z}_m$. The transposition class of the original pitch-class segment $X$ plays a privileged role in the definition of $\overline{Q}_j$, as we have $\overline{Q}_j(X) = T_j(X)$. For all other elements in $S$ we define $\overline{Q}_j$ contextually with respect to the linear parameter $a$ in $Y=( y_1,\dots, y_n ) = {^ba}(X) = T_b M_a (X) \in S$
$$\overline{Q}_j((y_1, y_2, \dots, y_m)) := T_{a j} ((y_1, y_2, \dots, y_m)).$$
\noindent Any of the Schritt Transformations commutes with any affine transformation  $f: \mathbb{Z}_{m} \to \mathbb{Z}_{m}$, say $f(y) = c y + d \,\, \mbox{mod} \,\, m$:

In order to contextually determine the image $\overline{Q}_j(f(Y))$ we first need to express the pitch-class segment $f(Y)$ as an image of $X$.  For all $i = 0, \dots, m-1$ one obtains:
$$f(y_i) = c y_i + d_i = c (a x_i + b) + d_i = (c a) x_i + c b + d_i.$$
In other words, the contextual multiplication factor is $a c$:
$$\begin{array}{lll}\overline{Q}_j(f(y_i)) & = & \overline{Q}_j(c y_i + d) \\
& = & T_{a c  j} (c y_i + d) \\
& = & c y_i + d + a c  j \\
& = & c (y_i + a j) + d \\
& = & f( \overline{Q}_j(y_i))
\end{array}$$

If $\mathrm{Aff}^\ast(\mathbb{Z}_{m})$ acts simply transitively on its orbit $S$ of
$X = (x_1,\dots, x_n)$, the contextual dual is $\{\overline{Q}_j W_a  \, | \, a \in \mathbb{Z}_m^\ast, j \in \mathbb{Z}_m \}.$ If we take $m=12$ and $X=(0,4,7)$, then this pair of 48-element dual groups contains the $T/I$-group and $PLR$-group as sub dual groups, as well as all three dualities listed in Figure~\ref{fig:dual_groups_to_G0}.

\section{Octatonic Example} \label{sec:Octatonic_Example}

We may now illustrate in our octatonic example the foregoing discussion of sub dual groups, generalized contextual groups, and commutativity of contextual inversions with affine maps.

Let $G$ be the mod 12 $T/I$-group.\footnote{In \cite{fiorenoll2011}, we used $H$ to denote the mod 12 $T/I$-group, rather than $G$. We use $G$ in this section in order to make the application of Theorem~\ref{thm:sub_dual_groups} apparent.} Consider the octatonic set $$O_{01}=\{0,1,3,4,6,7,9,10\}$$ in $\mathbb{Z}_{12}$ and its set-wise stabilizer $G_0$ in $G$, namely $$G_0=\{T_0,T_3,T_6,T_9,I_7,I_{10},I_1,I_4\}.$$
To any ordered tuple $(x_1,x_2,x_3)$ of three distinct pitch classes in the octatonic set, we have the associated contextual group $H^{(x_1,x_2,x_3)}$ acting on the 24-element set $$S^{(x_1,x_2,x_3)}:=T/I\text{-orbit of }(x_1,x_2,x_3)$$ as described in \cite{fioresatyendra2005} and recalled in Section~\ref{sec:Sub_Dual_Groups_and_Contextual_Groups}. Each contextual group $H^{(x_1,x_2,x_3)}$ is dihedral of order 24, and is dual to the mod 12 $T/I$-group $G$ in $\text{Sym}(S^{(x_1,x_2,x_3)})$. Note that we are considering the group $G$ embedded into $\text{Sym}(S^{(x_1,x_2,x_3)})$, so in this sense we speak about many groups that are dual to the $T/I$-group simultaneously, one for each $(x_1,x_2,x_3)$. We may now apply Theorem~\ref{thm:sub_dual_groups} to this situation.

We focus our attention on 3-tuples $(x_1,x_2,x_3)$ which are multiples of the major chord $\{0,4,7\}$ that are triads (as opposed to dyads). We list them and propose names in Figure~\ref{fig:047_multiples}. The major $(0,4,7)$, the jet $(0,4,1)$, and the stride $(0,4,10)$ are contained in the octatonic set $O_{01}$, therefore their orbits under the octatonic set stabilizer $G_0$ are also contained in $O_{01}$, as displayed in Figure~\ref{fig:dual_groups_to_G0}. The strain $(0,8,2)$, the shark $(0,8,11)$, and the minor $(0,8,5)$ are not in $O_{01}$, but their translations by 1 are.
\begin{figure}
\caption{Multiples of $\{0,4,7\}$ which have three notes} \label{fig:047_multiples}
$$
\begin{array}{lcll}
\{0,4,7\}\times 1 &=&\{0,4,7\} &=\text{ major type}\\
\{0,4,7\}\times 2 &=&\{0,8,2\} &=\text{ strain type}\\
\{0,4,7\}\times 5 &=&\{0, 8, 11\} &=\text{ shark type}\\
\{0,4,7\}\times 7 &=&\{0,4,1\}&=\text{ jet type }\\
\{0,4,7\}\times 10 &=&\{0,4,10\}&=\text{ stride type}\\
\{0,4,7\}\times 11 &=&\{0,8,5\}&=\text{ minor type}
\end{array}
$$
\end{figure}

\begin{figure}
\begin{center}
\begin{tabular}{|c|l|c|}
\hline
\multicolumn{3}{|c|}{Simply Transitive Covers of Octatonic $O_{01}=\{0,1,3,4,6,7,9,10\}$ } \\
\hline
$s_0=(x_1,x_2,x_3)$ & $S_0^{(x_1,x_2,x_3)}=G_0\text{-orbit of }(x_1,x_2,x_3)$ & Dual Group \\ \hline \hline
$(0,4,7)$ & Major: (0, 4, 7), (3, 7, 10), (6, 10, 1), (9, 1, 4), & $H^{(0,4,7)}_0$ \\
& Minor: (7, 3, 0), (10, 6, 3), (1, 9, 6), (4, 0, 9) & \\ \hline
$(0,4,1)$ & Jet:\phantom{ark} (0, 4, 1), (9, 1, 10), (6, 10, 7), (3, 7, 4)  & $H^{(0,4,1)}_0$ \\
 & Shark: (1, 9, 0), (10, 6, 9), (7, 3, 6), (4, 0, 3) & \\ \hline
$(0, 4, 10)$ & Stride: (0, 4, 10), (3, 7, 1), (6, 10, 4), (9, 1, 7) & $H^{(0,4,10)}_0$ \\
& Strain: (10, 6, 0), (1, 9, 3), (4, 0, 6), (7, 3, 9) &  \\ \hline
\end{tabular}
\end{center}
\caption{Dual groups to $G_0=\{T_0,T_3,T_6,T_9,I_7,I_{10},I_1,I_4\}$. See also Figure~\ref{fig:047_M7_M10_network} for a diagrammatic presentation.} \label{fig:dual_groups_to_G0}
\end{figure}

We next consider the group structure of $G_0$ and its dual groups $H_0^{(0,4,7)}$, $H_0^{(0,4,1)}$, and $H_0^{(0,4,10)}$. The group $G_0$ is dihedral of order 8, since the dihedral relations hold: $s:=T_3$ has order 4, $t:=I_1$ has order 2, and $tst=s^{-1}$ holds. The three dual groups $H_0^{(0,4,7)}$, $H_0^{(0,4,1)}$, and $H_0^{(0,4,10)}$ are also dihedral of order 8, as they are isomorphic to $G_0$. To determine the elements of the dual groups,  by Theorem~\ref{thm:sub_dual_groups}~\ref{thm:sub_dual_groups:H_0_simple_transitivity} and \ref{thm:sub_dual_groups:restricted_subgroups_dual}, we must only determine which elements of $H^{(x_1,x_2,x_3)}$ map the first column in Figure~\ref{fig:dual_groups_to_G0} to the respective elements in the second column of Figure~\ref{fig:dual_groups_to_G0}. We see that the various powers and composites of $Q_3$ and $P$-analogues reach all eight elements, so the dual groups are (the restrictions of)
$$
\aligned
H_0^{(0,4,7)} &= \{Q_0,Q_3,Q_6,Q_9,P^{(0,4,7)},Q_3P^{(0,4,7)},Q_6P^{(0,4,7)},Q_9P^{(0,4,7)} \} \\
H_0^{(0,4,1)} &= \{Q_0,Q_3,Q_6,Q_9,P^{(0,4,1)},Q_3P^{(0,4,1)},Q_6P^{(0,4,1)},Q_9P^{(0,4,1)}\} \\
H_0^{(0,4,10)} &= \{Q_0,Q_3,Q_6,Q_9,P^{(0,4,10)},Q_3P^{(0,4,10)},Q_6P^{(0,4,10)},Q_9P^{(0,4,10)}\}. \\
\endaligned
$$
The notation $P^{(x_1,x_2,x_3)}$ indicates the $P$-analogue for the pitch-class segment $(x_1,x_2,x_3)$, which is $J^{1,3}$ (see equation \eqref{equ:contextual_inversion} and the discussion thereafter). We follow the same convention $R^{(x_1,x_2,x_3)}=J^{1,2}$ below. This definition
of $P$ and $R$ on the $T/I$-orbits of the pitch-class segments ordered as suggested in Figure~\ref{fig:047_multiples} is reasonable because they directly correlate with $P$ and $R$ on the $T/I$-orbit of $(0,4,7)$.

From Figure~\ref{fig:dual_groups_to_G0} and Theorem~\ref{thm:sub_dual_groups}~\ref{thm:sub_dual_groups:H_0_simple_transitivity} we also see that the $R$-analogues must be in the dual groups, so $R$ can be written in terms of $Q_3$ and $P$. Computing\footnote{To compute the composite $PR$, one can use the explicit formulas for $P$ and $R$, or one could simply evaluate on $(0,4,7)$, $(0,4,1)$, $(0,4,10)$ and use the uniqueness guaranteed by the simple transitivity of the respective $H$-actions.}, we see
$$\aligned
P^{(0,4,7)}R^{(0,4,7)}&=Q_9 \\
P^{(0,4,1)}R^{(0,4,1)}&=Q_3 \\
P^{(0,4,10)}R^{(0,4,10)}&=Q_6. \\
\endaligned$$
Therefore, $\langle P^{(0,4,7)},R^{(0,4,7)} \rangle=H_0^{(0,4,7)}$ and $\langle P^{(0,4,1)},R^{(0,4,1)} \rangle=H_0^{(0,4,1)}$, while $\langle P^{(0,4,10)},R^{(0,4,10)} \rangle$ is a Klein 4-group properly contained in $H_0^{(0,4,10)}$.

In Figures~\ref{fig:P-analogues} and \ref{fig:R-analogues}, we write out the $P$- and $R$-analogues for the triads in Figure~\ref{fig:dual_groups_to_G0}, using $(a,b,c)$ in place of $(y_1,y_2,y_3)$ for readability.

\begin{figure}
$$\begin{array}{lll}
P^{(0, 4, 7)}(a, b, c) & := & (c, b-1, a) \quad \mbox{if } (a, b, c) \in Maj,\\
P^{(0, 4, 7)}(a, b, c) & := & (c, b+1, a) \quad \mbox{if } (a, b, c) \in Min, \\
\\
P^{(0, 4, 1)}(a, b, c) & := & (c, b+5, a) \quad \mbox{if } (a, b, c) \in Jet,\\
P^{(0, 4, 1)}(a, b, c) & := &(c, b-5, a) \quad \mbox{if } (a, b, c) \in Shark, \\
\\
P^{(0, 4, 10)}(a, b, c) & := & (c, b+2, a) \quad \mbox{if } (a, b, c) \in Stride,\\
P^{(0, 4, 10)}(a, b, c) & := & (c, b-2, a) \quad \mbox{if } (a, b, c) \in Strain, \\
\\
\end{array}$$
\caption{$P$-analogues for the triads in Figure~\ref{fig:dual_groups_to_G0}} \label{fig:P-analogues}
\end{figure}

\begin{figure}
$$\begin{array}{lll}
R^{(0, 4, 7)}(a, b, c) & := & (b, a, c +2) \quad \mbox{if } (a, b, c) \in Maj,\\
R^{(0, 4, 7)}(a, b, c) & := & (b, a, c -2) \quad \mbox{if } (a, b, c) \in Min, \\
\\
R^{(0, 4, 1)}(a, b, c) & := & (b, a, c +2) \quad \mbox{if } (a, b, c) \in Jet,\\
R^{(0, 4, 1)}(a, b, c) & := & (b, a, c -2) \quad \mbox{if } (a, b, c) \in Shark, \\
\\
R^{(0, 4, 10)}(a, b, c) & := & (b, a, c - 4) \quad \mbox{if } (a, b, c) \in Stride,\\
R^{(0, 4, 10)}(a, b, c) & := & (b, a, c + 4) \quad \mbox{if } (a, b, c) \in Strain, \\
\\
\end{array}$$
\caption{$R$-analogues for the triads in Figure~\ref{fig:dual_groups_to_G0}} \label{fig:R-analogues}
\end{figure}

Pictured in Figure~\ref{fig:dual_groups_to_T1G0T1Inverse} is the $T_1$-transform of the systems in Figure~\ref{fig:dual_groups_to_G0}.
Figure~\ref{fig:Figure1OpeningTheme} shows the opening theme of Schoenberg's op. 7 quartet, which, together with the affine image in the lower staff of Figure~\ref{fig:Figure7PiecewideNarrative}, engages the systems of Figure~\ref{fig:dual_groups_to_T1G0T1Inverse}.
Recall from \cite[Corollary 3.3]{fiorenoll2011} that if $k \in G$, then the orbit of $ks_0$ under $H_0$  is $kS_0$ and the associated dual group of $H_0$ is $kG_0k^{-1}$ (the roles of $H_0$ and $G_0$ are switched in this paraphrase because that is what the analysis requires).  The opening theme takes place in the octatonic set
$$O_{12}=\{1,2,4,5,7,8,10,11\}$$
and the sets $S_0^{(x_1,x_2,x_3)}$ of Figure~\ref{fig:dual_groups_to_T1G0T1Inverse} are simply transitive covers of $O_{12}$.
\begin{figure}
\begin{center}
\begin{tabular}{|c|l|c|}
\hline
\multicolumn{3}{|c|}{Simply Transitive Covers of Octatonic $O_{12}=\{1,2,4,5,7,8,10,11\}$ } \\
\hline
$s_0=(x_1,x_2,x_3)$ & $S_0^{(x_1,x_2,x_3)}=G_0\text{-orbit of }(x_1,x_2,x_3)$ & Dual Group \\ \hline \hline
$(1,5,8)$ & Major: (1, 5, 8), (4, 8, 11), (7, 11, 2), (10, 2, 5), & $H^{(1,5,8)}_0$ \\
& Minor: (8, 4, 1), (11, 7, 4), (2, 10, 7), (5, 1, 10) & \\ \hline
$(1, 5, 2)$ & Jet:\phantom{ark} (1, 5, 2), (10, 2, 11), (7, 11, 8), (4, 8, 5)  & $H^{(1,5,2)}_0$ \\
 & Shark: (2, 10, 1), (11, 7, 10), (8, 4, 7), (5, 1, 4) & \\ \hline
$(1, 5, 11)$ & Stride: (1, 5, 11), (4, 8, 2), (7, 11, 5), (10, 2, 8) & $H^{(1,5,11)}_0$ \\
& Strain: (11, 7, 1), (2, 10, 4), (5, 1, 7), (8, 4, 10) &  \\ \hline
\end{tabular}
\end{center}
\caption{Dual groups to $T_1G_0T_1^{-1}=\{T_0,T_3,T_6,T_9,I_9,I_{0},I_3,I_6\}$} \label{fig:dual_groups_to_T1G0T1Inverse}
\end{figure}
The operations in the groups $H_0^{(1,5,8)}$, $H_0^{(1,5,2)}$, and $H_0^{(1,5,11)}$ of Figure~\ref{fig:dual_groups_to_T1G0T1Inverse} are the restrictions of the same operations as those in the respective groups $H_0^{(0,4,7)}$, $H_0^{(0,4,1)}$, and $H_0^{(0,4,10)}$ displayed above (they are not conjugated, only $G_0$ is conjugated). Consequently, the dual groups $H_0^{(1,5,8)}$, $H_0^{(1,5,2)}$, and $H_0^{(1,5,11)}$ in Figure~\ref{fig:dual_groups_to_T1G0T1Inverse}, are also each generated by the respective $P$ and $R$ analogues for their pitch-class segments.

Returning to the $O_{01}$ systems of Figure~\ref{fig:dual_groups_to_G0}, we have morphisms $M_7$ and $M_{10}$ from the major/minor dual pair to the jet/shark and stride/strain dual pairs.  By Corollary~\ref{cor:affines_commute_with_subgroups_generated_by_contextual_inversions}, multiplication by 7 and 10 induce morphisms of the larger simply transitive group actions
$$\xymatrix{(M_7,\varphi)\colon (H^{(0,4,7)},S^{(0,4,7)}) \ar[r] & (H^{(0,4,1)},S^{(0,4,1)})}$$
$$\xymatrix{(M_{10},\psi)\colon (H^{(0,4,7)},S^{(0,4,7)}) \ar[r] & (H^{(0,4,10)},S^{(0,4,10)})},$$
because $H^{(0,4,7)}=\langle L,R \rangle$.
These also restrict to the sub simply transitive group actions
$$\xymatrix{(M_7,\varphi)\colon(H^{(0,4,7)}_0,S^{(0,4,7)}_0) \ar[r] & (H^{(0,4,1)}_0,S^{(0,4,1)}_0)}$$
$$\xymatrix{(M_{10},\psi)\colon (H^{(0,4,7)}_0,S^{(0,4,7)}_0) \ar[r] & (H^{(0,4,10)}_0,S^{(0,4,10)}_0)}.$$
and induce morphisms of the corresponding generalized interval systems. The morphism $(M_{10},\psi)$ and its restriction are not epic, as $M_{10}$ is not surjective (coordinates in its image are always even), so $\psi$ cannot be surjective by the simply transitive counterpart to Proposition~\ref{prop:monic/epic_morphisms_for_intervallic}. Since $S^{(0,4,7)}$ and $S^{(0,4,1)}$ have the same finite cardinality, namely 24, $M_{10}$ cannot be injective either. Note also $L^{(0,4,10)}R^{(0,4,10)}=Q_2$ and $\langle L^{(0,4,10)},R^{(0,4,10)} \rangle$ is properly contained in $H^{(0,4,10)}$. In other words, the generalized contextual group for $(0,4,10)$ is not generated by $L$ and $R$.

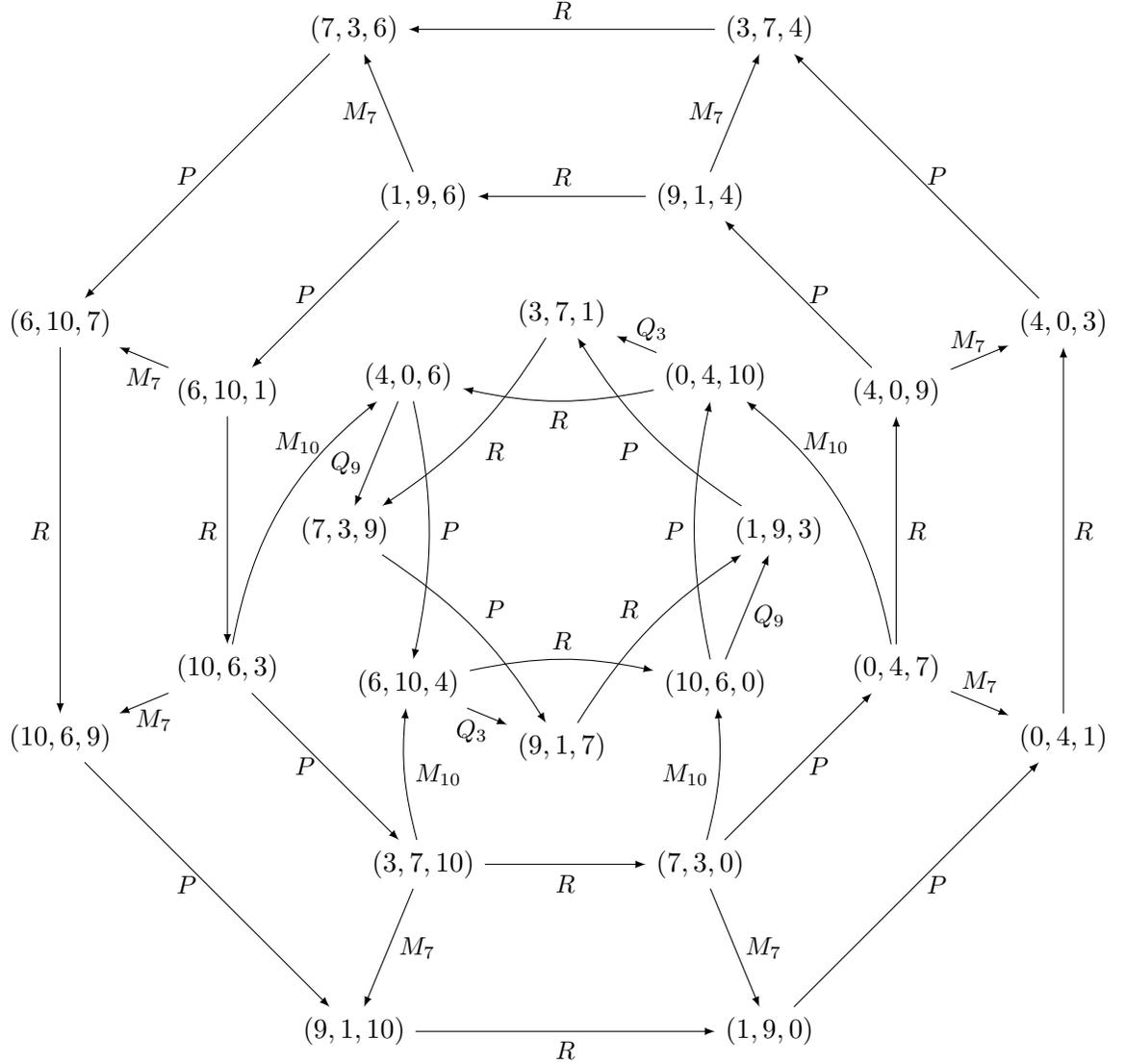
\begin{figure}
\caption{$PR$-Networks for the major $(0,4,7)$ and the jet $(0,4,1)$, and a $PRQ$-network for the stride $(0,4,10)$, all from Figure~\ref{fig:dual_groups_to_G0}, connected by GIS morphisms $M_7$ and $M_{10}$ } \label{fig:047_M7_M10_network}

\begin{flushleft}
\begin{tikzpicture}[>=latex]

  \def\radius{7.5 cm} 
    \node[font=\small] (104) at (-22.5:\radius)   {$(0, 4, 1)$};
    \node[font=\small]  (043) at (22.5:\radius)  {$(4, 0, 3)$};
    \node[font=\small]  (437) at (67.5:\radius) {$(3, 7, 4)$};
    \node[font=\small]  (376) at (112.5:\radius)  {$(7, 3, 6)$};
    \node[font=\small]  (76t) at (157.5:\radius)  {$(6, 10, 7)$};
    \node[font=\small]  (6t9) at (202.5:\radius)  {$(10, 6, 9)$};
    \node[font=\small]  (t91) at (247.5:\radius)  {$(9, 1, 10)$};
    \node[font=\small]  (910) at (292.5:\radius)  {$(1, 9, 0)$};

 \def\radius{5.0 cm} 
    \node[font=\small]  (704) at (-22.5:\radius)   {$(0, 4, 7)$};
    \node[font=\small]  (049) at (22.5:\radius)  {$(4, 0, 9)$};
    \node[font=\small]  (491) at (67.5:\radius) {$(9, 1, 4)$};
    \node[font=\small]  (916) at (112.5:\radius)  {$(1, 9, 6)$};
    \node[font=\small]  (16t) at (157.5:\radius)  {$(6, 10, 1)$};
    \node[font=\small]  (6t3) at (202.5:\radius)  {$(10, 6, 3)$};
    \node[font=\small]  (t37) at (247.5:\radius)  {$(3, 7, 10)$};
    \node[font=\small]  (370) at (292.5:\radius)  {$(7, 3, 0)$};

\def\radius{3 cm} 
    \node[font=\small]  (913) at (0:\radius)   {$(1, 9, 3)$};
    \node[font=\small]  (137) at (90:\radius)  {$(3, 7, 1)$};
    \node[font=\small]  (379) at (180:\radius) {$(7, 3, 9)$};
    \node[font=\small]  (791) at (270:\radius)  {$(9, 1, 7)$};

 \def\radius{3 cm} 
    \node[font=\small]  (t04) at (45:\radius)  {$(0, 4, 10)$};
    \node[font=\small]  (046) at (135:\radius)  {$(4, 0, 6)$};
    \node[font=\small]  (46t) at (225:\radius)  {$(6, 10, 4)$};
    \node[font=\small]  (6t0) at (315:\radius)  {$(10, 6, 0)$};

 \path[->,font=\footnotesize]

       (704) edge node[above, pos=.5] {$M_7$} (104)
        (049) edge node[above, pos=.3] {$M_7 $} (043)
        (491) edge  node[left] {$M_7 $} (437)
        (916) edge node[left, pos=.5] {$M_7 $} (376)
        (16t) edge node[below, pos=.5] {$M_7 $} (76t)
        (6t3) edge node[below, pos=.3] {$M_7 $} (6t9)
        (t37) edge node[right] {$M_7 $} (t91)
        (370) edge node[right, pos=.5] {$M_7 $} (910)

        (704) edge [bend right = 20] node[right, pos=.8] {$M_{10}$} (t04)
        (6t3) edge [bend left = 20] node[left, pos=.8] {$M_{10}$} (046)
        (t37) edge [bend left = 10] node[right] {$M_{10}$} (46t)
        (370) edge [bend right = 10] node[left] {$M_{10}$} (6t0)

         (t04) edge [bend right = 0] node[above, pos=.1] {$Q_3$} (137)
        (046) edge  [bend right = 0] node[left, pos=.6] {$Q_9$} (379)
        (46t) edge [bend right = 0] node[below, pos=.1] {$Q_3$} (791)
       (6t0) edge [bend left = 0] node[right, pos=.4] {$Q_9$} (913)

         (104) edge node[right] {$R$} (043)
        (043) edge node[right] {$P$} (437)
        (437) edge node[above] {$R$} (376)
        (376) edge node[left] {$P$} (76t)
        (76t) edge  node[left] {$R$} (6t9)
        (6t9) edge node[left] {$P$} (t91)
        (t91) edge node[below] {$R$} (910)
        (910) edge node[right] {$P$} (104)

        (704) edge node[right] {$R$} (049)
        (049) edge node[right] {$P$} (491)
        (491) edge  node[above] {$R$} (916)
        (916) edge node[left] {$P$} (16t)
        (16t) edge node[left] {$R$} (6t3)
        (6t3) edge node[left] {$P$} (t37)
        (t37) edge node[below] {$R$} (370)
        (370) edge node[right] {$P$} (704)

        (137) edge [bend left = 12]  node[auto] {$R$} (379)
        (379) edge [bend left = 12]  node[auto] {$P$} (791)
        (791) edge [bend left = 12]  node[auto] {$R$} (913)
        (913) edge [bend left = 12]  node[auto] {$P$} (137)
        (t04) edge   [bend left = 12]  node[auto] {$R$} (046)
        (046) edge  [bend left = 12]   node[auto] {$P$} (46t)
        (46t) edge  [bend left = 12]  node[auto] {$R$} (6t0)
        (6t0) edge   [bend left = 12]  node[auto] {$P$} (t04);
\end{tikzpicture}
\end{flushleft}

\end{figure}

Figure~\ref{fig:047_M7_M10_network} shows the simply transitive group actions in the octatonic $O_{01}$ of Figure~\ref{fig:dual_groups_to_G0}, connected by the morphisms $M_7$ and $M_{10}$. The two octagons illustrate that the $PR$-groups for the major $(0,4,7)$ and the jet $(0,4,1)$ are dihedral of order 8, while the two squares in the middle illustrate that the $PR$-group for the stride $(0,4,10)$ is a Klein 4-group. The $(0,4,10)$ group elements $Q_3$ and $Q_9$ connect the two squares to make a cube, as guaranteed by simple transitivity.  This network, on the one hand, paradigmatically saturates the analytical Summary Network from Figure~\ref{fig:Figure8SchoenbergNetworkRP} and, on the other hand, relocates it within a single octatonic collection.  The top row of the Summary Network is the outer octagon transposed by 1, the second and fourth rows are the inner octagon, while the two bottom rows are part of the inner squares. The Summary Network is a further refinement of Figure~\ref{fig:047_M7_M10_network} in that the orderings of the pitch-class segments correspond exactly to the temporal order in Schoenberg's String Quartet in $D$ minor.

The morphism ${}^77=T_7M_7$ is central to our analysis of Schoenberg's quartet. It is a morphism from the dual major/minor subsystem in Figure~\ref{fig:dual_groups_to_G0} to the dual shark/jet subsystem in Figure~\ref{fig:dual_groups_to_T1G0T1Inverse}
$$\xymatrix{{}^77\co (H_0^{(0,4,7)},S_0^{(0,4,7)}) \ar[r] & (H_0^{(1,5,2)},S_0^{(1,5,2)})}$$
because ${}^77(0,4,7)=(7,11,8)\in S_0^{(1,5,2)}$, both $H_0^{(0,4,7)}$ and  $H_0^{(1,5,2)}$ are generated by $P$ and $R$ analogues, and Corollary~\ref{cor:affines_commute_with_subgroups_generated_by_contextual_inversions}
applies. Note that the pitch-class sets $\{0,4,7\}$ and $\{1,5,2\}$ are asymmetric, so ``forgetting order'' defines bijections
$$\xymatrix{\text{$T/I$-class of }(0,4,7) \ar[r] & \text{$T/I$-class of } \{0,4,7\}}$$
$$\xymatrix{\text{$T/I$-class of }(1,5,2) \ar[r] & \text{$T/I$-class of } \{1,5,2\}},$$
which induce actions of $H_0^{(0,4,7)}$ and  $H_0^{(1,5,2)}$ on the respective subcollections of {\it unordered} pitch-class sets, and we have ${}^77$ as a morphism of these unordered systems as well.

However, in Section~\ref{sec:Analysis_of_Schoenberg} we will not need the unordered systems. Instead, we use permutations in a more refined simply transitive action so that the order in the pitch-class segments corresponds to the local temporal order in the work of music. Permutations in neo-Riemannian theory were developed by Fiore--Noll--Satyendra in \cite{fiorenollsatyendraMCM2013}, and are sketched at the end of Section~\ref{sec:Analysis_of_Schoenberg} to make the present paper self-contained.

\section{Elements of a Transformational Analysis of Schoenberg, String Quartet in $D$ Minor, Op. 7} \label{sec:Analysis_of_Schoenberg}

\begin{figure}[h]\caption{Schoenberg, String Quartet in $D$ minor, op. 7. Opening theme materials, (measures 1-3, violin 1)} \label{fig:Figure1OpeningTheme}
\includegraphics[height=1.25in]{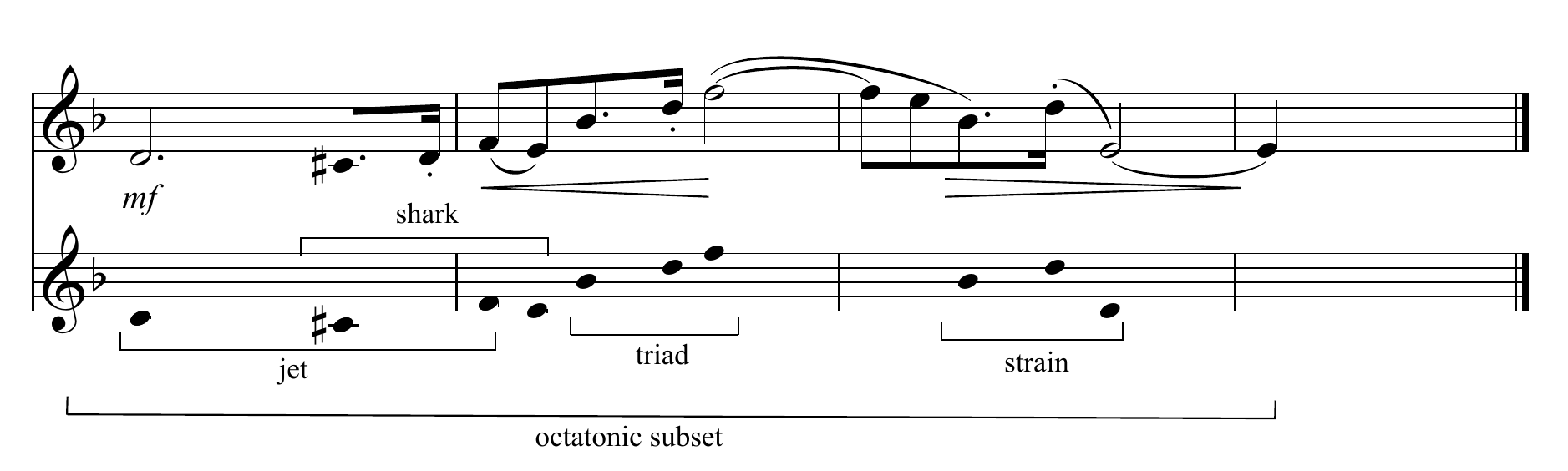}
\end{figure}

In this section we illustrate the analytical potential of the results of Sections~\ref{sec:Morphisms}, \ref{sec:Sub_Dual_Groups_and_Contextual_Groups}, and \ref{sec:Octatonic_Example} by integrating selected observations about the motivic organization of Schoenberg's String Quartet in $D$ minor, op. 7 into an analytical network. In particular, we use morphisms to indicate a parallel organization between a shark-jet cycle and a major-minor triadic cycle, in which each cycle is an octatonic cover.
The four set-types jet, shark, major, and strain in the $O_{12}$ octatonic set in Figure~\ref{fig:dual_groups_to_T1G0T1Inverse} are melodically expressed within the opening violin melody, as indicated in Figure~\ref{fig:Figure1OpeningTheme}.

\begin{figure}\caption{Schoenberg, String Quartet in $D$ minor, op. 7.  Relations between shark, jet, strain, and triad forms in the opening theme} \label{fig:Figure2OpeningNetwork}
$$
\renewcommand{\labelstyle}{\textstyle}
\entrymodifiers={+[F-]}
\xymatrix@C=4pc@R=3pc{(2,1,5) \ar[r]^{\mathrm{RICH}}_{(13)R} & (1,5,4) \ar[r]^{{}^9 5} & (2,10,5) \ar[r]^{{}^6 2} & (10,2,4)}$$
\end{figure}
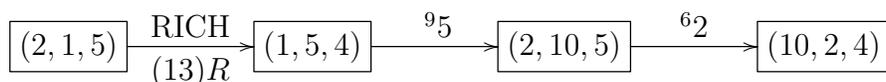

The operations on the ordered pitch-class sets in Section~\ref{sec:Octatonic_Example} induce operations on their reordered variants in this analysis, via permutations. For instance we may say we have an {\it $R$-relation} between the jet $(2,1,5)$ and the shark $(1,5,4)$ in the opening theme because $R(1,5,2)=(5,1,4)$ in the ordered context of the $(H_0^{(1,5,2)}, S_0^{(1,5,2)})$-system of Figure~\ref{fig:dual_groups_to_T1G0T1Inverse} in Section~\ref{sec:Octatonic_Example} (for readability we do not write the superscripts on $R$ anymore).  Intuitively, the $R$-operation on a ``non-root position'' trichord is defined by permuting to ``root position'', applying $R$, and then permuting back.  Notice however that the entries of the jet $(2,1,5)$ and the shark $(1,5,4)$, which appear exactly in the score in this order, do not correspond to the entries of $R(1,5,2)=(5,1,4)$ in the same way: the first entry of $(2,1,5)$ is the last entry of $(1,5,2)$, whereas the first entry of $(1,5,4)$ is not the last entry of $(5,1,4)$. Thus, the transformation from $(2,1,5)$ to $(1,5,4)$ is more than a mere (conjugate of the) $R$ transformation; it in fact also has a permutational effect which interchanges the first and the third components. We have indicated this in Figure~\ref{fig:Figure2OpeningNetwork} as $(13)R$ where $(13)$ is the cycle notation for the permutation on three ``letters'' which switches 1 and 3 and fixes 2. For more details, see the end of this section and also \cite{fiorenollsatyendraMCM2013}.

Returning now to the analysis, Figure~\ref{fig:Figure2OpeningNetwork} labels transformational relations between the adjacent set types of the opening theme; the arrow direction denotes temporal order. Recall that $\mathrm{RICH}$ in Figure~\ref{fig:Figure2OpeningNetwork} means {\it retrograde inversion enchaining},\footnote{We use the term ``retrograde'' in a more general sense of ``reversal''.} that is, if $Y$ is a pitch-class segment, then $\mathrm{RICH}(Y)$ is that retrograde inversion of $Y$ which has the first two notes $y_{n-1}$ and $y_{n}$, in that order. $\mathrm{RICH}$ is not part of the groups of Section~\ref{sec:Octatonic_Example}, see \cite{fiorenollsatyendraMCM2013} for a clarification of the relationship between $\mathrm{RICH}$ and permutations. For the first two chords in Figure~\ref{fig:Figure2OpeningNetwork}, $\mathrm{RICH}$ is an $R$-relation.  As an unrelated remark, we have written $(2,10,5)$ instead of the triad $(10,2,5)$ (as in the music) to illustrate the componentwise correspondence under the operation ${}^95$.

\begin{figure}\caption{Schoenberg, String Quartet in $D$ minor, op. 7.  Reduction of measures 8-10} \label{fig:Figure3ReductionOf8-10ToRP}
\includegraphics[height=1.25in]{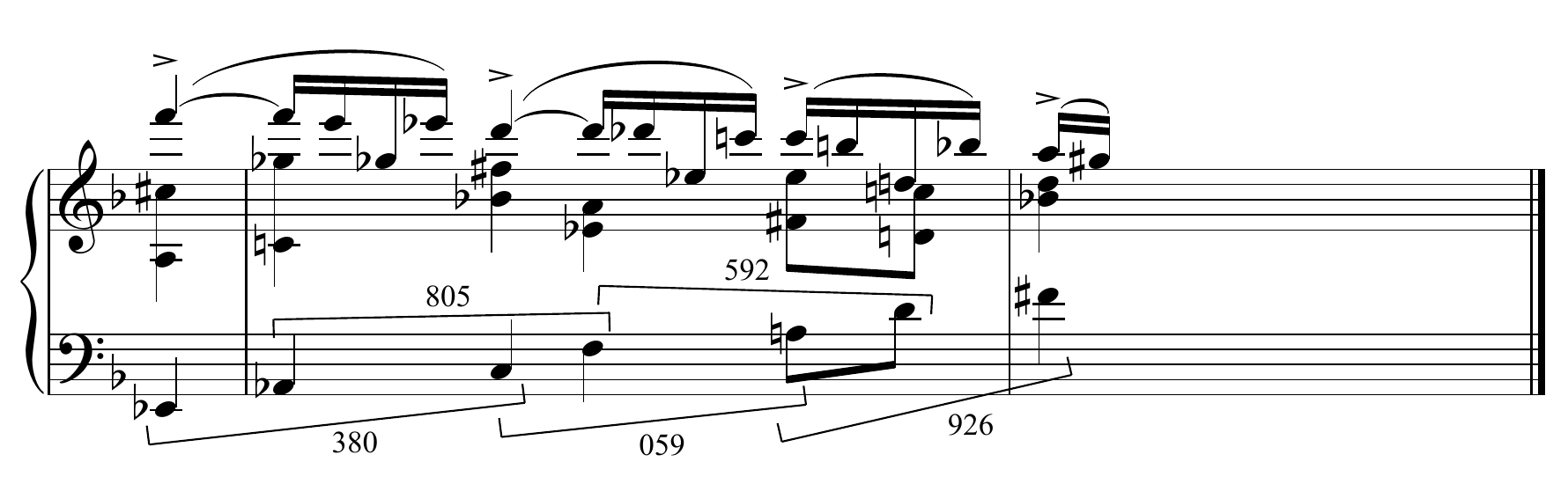}
\end{figure}

\begin{figure}\caption{$\mathrm{RICH}$-chaining on the triadic motive in measures 8-10, cello} \label{fig:Figure4RICHChaining8-10}
$$\renewcommand{\labelstyle}{\textstyle}
\entrymodifiers={+[F-]}
\xymatrix@C=4pc@R=3pc{(3,8,0) \ar[r]^{\mathrm{RICH}}_{(13)R} & (8,0,5) \ar[r]^{\mathrm{RICH}}_{(13)P} & (0,5,9) \ar[r]^{\mathrm{RICH}}_{(13)R} & (5,9,2) \ar[r]^{\mathrm{RICH}}_{(13)P} & (9,2,6) }$$
\end{figure}
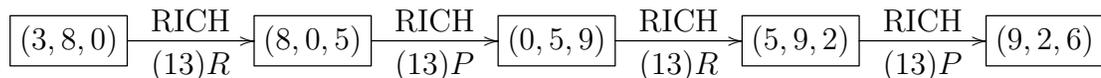

Figure~\ref{fig:Figure3ReductionOf8-10ToRP} is a reduction of an octatonic passage that dramatically interrupts the $D$-minor diatonicism of the opening.  The dominants on $A\flat$, $F$, and $D$ lend a strong octatonic flavor to the harmony, which is not disturbed by the extraneous $C\sharp$ and $B\flat$ in the middle voices. (As we will see later, these extraneous notes create whole-tone harmonies that are thematized in transformations between stride and strain trichords.)
Of particular significance for our purposes is the octatonic cello melody in the lower staff of Figure~\ref{fig:Figure3ReductionOf8-10ToRP} that is analyzed in Figure~\ref{fig:Figure4RICHChaining8-10} as a $\mathrm{RICH}$-chain.  As shown, the melody may also be analyzed as an octatonic $RP$-chain of neo-Riemannian operations, see Cohn's article \cite{cohn1997}. Figure~\ref{fig:Figure4RICHChaining8-10} takes place on the ocatonic $O_{23}$ in the simply transitive action $(H_0^{(2,6,9)}, S_0^{(2,6,9)})$ , which is a transform of $(H_0^{(0,4,7)}, S_0^{(0,4,7)})$ by $T_2$.

In measures 8-10 (see Figures~\ref{fig:Figure3ReductionOf8-10ToRP} and \ref{fig:Figure4RICHChaining8-10}), the $\mathrm{RICH}$-chaining does not complete the eight-triad cycle; the cello melody breaks off after five triads. A complete cycle of eight triads is traversed in the octatonic $O_{01}$ by the simply transitive action $(H_0^{(0,4,7)},S_0^{(0,4,7)})$, however, in an extraordinary passage that ends the first large formal section of the quartet (measures 1-96):  $G\flat$, $E\flat$-minor, $E\flat$, $C$-minor, $C$, $A$-minor, $A$, $F\sharp$-minor.  Figure~\ref{fig:Figure5CompleteRPChain88-92} shows the triadic melody along with a harmonic reduction of chords that support it, which can be analyzed by register as expressing $\mathrm{RICH}$-related stride-strain pairs: $(7,3,1)$-$(3,1,9)$ and $(4,0,10)$-$(0,10,6)$.

Figure~\ref{fig:Figure6TransformationOfPRCycleNotes} shows an equivalence between the shark-jet and major-minor motivic organization in the quartet. Applying the morphism
$$\xymatrix{{}^77\co (H_0^{(0,4,7)},S_0^{(0,4,7)}) \ar[r] & (H_0^{(1,5,2)},S_0^{(1,5,2)})}$$
to the major-minor series in the upper staff gives the shark-jet series in the lower staff.

The unfolding of the shark-jet cycle may be heard as a piece-wide narrative, as suggested by the following discussion of Figure~\ref{fig:Figure7PiecewideNarrative}. The opening violin motive (Figures~\ref{fig:Figure1OpeningTheme} and \ref{fig:Figure7PiecewideNarrative}), gives forms 1 and 2. At C30, the $\mathrm{RICH}$-chain continues with forms 3 and 4 in a {\it fortissimo} return of the main theme doubled in all four instruments.   At measure 85 just before the complete triadic cycle--an important moment as it ends the first formal unit spanning measures 1-94--forms 5, 6, and 7 are highlighted through imitation between violins 1 and 2. An extraordinary return of the $D$-minor theme in the parallel major concludes the work.  In these final bars the final forms, 7 and 8, are articulated as the theme's registral extremes.  The work concludes shortly thereafter.

\begin{figure}\caption{Summary Network, Schoenberg, String Quartet in $D$ minor, op. 7. Every horizontal arrow could also be labelled $\mathrm{RICH}$. Upper labels on boxes are measure numbers.} \label{fig:Figure8SchoenbergNetworkRP}
{\footnotesize
$$
\renewcommand{\labelstyle}{\textstyle}
\entrymodifiers={+[F-]}
\xymatrix@C=1pc@R=3pc{ *\txt{shark-jet \\ motive} & (2,1,5) \ar[r]_*!/_7pt/{(13)R} {}\save[]+<-0.4cm,.5cm>*\txt<8pc>{%
1--2} \restore & (1,5,4) \ar[r]_*!/_7pt/{(13)P} & (5,4,8) \ar[r]_*!/_7pt/{(13)R} {}\save[]+<-.4cm,.5cm>*\txt<8pc>{%
C30} \restore & (4,8,7) \ar[r]_*!/_7pt/{(13)P} & (8,7,11) \ar[r]_*!/_7pt/{(13)R} {}\save[]+<-0.6cm,.5cm>*\txt<8pc>{%
85} \restore & (7,11,10) \ar[r]_*!/_7pt/{(13)P} & (11,10,2) \ar[r]_*!/_7pt/{(13)R} {}\save[]+<-0.2cm,.5cm>*\txt<8pc>{%
O13--15} \restore  & (10,2,1) \\
*\txt{triadic \\ motive \\ v1, vc, v1} & (1,6,10) \ar[r]_*!/_7pt/{(13)R} \ar[u]_{{}^77} \ar[d]^{{}^21} {}\save[]+<-0.5cm,.5cm>*\txt<8pc>{%
88--92} \restore & (6,10,3) \ar[r]_*!/_6pt/{(13)P} \ar[u]_{{}^77} \ar[d]^{{}^21} & (10,3,7) \ar[r]_*!/_7pt/{(13)R} \ar[u]_{{}^77} \ar[d]^{{}^21} & (3,7,0) \ar[r]_*!/_7pt/{(13)P} \ar[u]_{{}^77} \ar[d]^{{}^21} & (7,0,4) \ar[r]_*!/_7pt/{(13)R} \ar[u]_{{}^77} \ar[d]^{{}^21} & (0,4,9) \ar[r]_*!/_7pt/{(13)P} \ar[u]_{{}^77}  & (4,9,1) \ar[r]_*!/_7pt/{(13)R} \ar[u]_{{}^77} & (9,1,6) \ar[u]_{{}^77} \\
*\txt{triadic \\ motive \\ vc} & (3,8,0) \ar[r]_*!/_6pt/{(13)R} \ar[d]^{ {}^41} {}\save[]+<-0.5cm,.5cm>*\txt<8pc>{%
8--10} \restore & (8,0,5) \ar[r]_*!/_7pt/{(13)P} \ar[d]^{ {}^41} & (0,5,9) \ar[r]_*!/_7pt/{(13)R} \ar[d]^{ {}^41} & (5,9,2) \ar[r]_*!/_7pt/{(13)P} \ar[d]^{ {}^41} & (9,2,6) \ar[d]^{ {}^41} & *{} & *{} & *{} \\
 *\txt{triadic \\ motive \\ vc} & (7,0,4) \ar[r]_*!/_7pt/{(13)R} \ar[d]^{{}^010} {}\save[]+<-0.6cm,.5cm>*\txt<8pc>{%
10--12} \restore & (0,4,9) \ar[r]_*!/_7pt/{(13)P} \ar[d]^{{}^010} & (4,9,1) \ar[r]_*!/_7pt/{(13)R} \ar[d]^{{}^010} & (9,1,6) \ar[r]_*!/_7pt/{(13)P}
& (1,6,10) & *{} & *{} & *{} \\
*\txt{stride-strain \\ motive \\ v2-vc} & (10,0,4) \ar[d]^{{}^31} {}\save[]+<-0.6cm,.5cm>*\txt<8pc>{%
91} \restore & (0,4,6) \ar[l]^*!/^7pt/{(13)R} \ar[d]^{{}^31} & (4,6,10) \ar[l]^*!/^7pt/{(13)P} \ar[d]^{{}^31} \\
*\txt{stride-strain \\ motive \\ v1-va} & (1,3,7) \ar[r]_*!/_7pt/{(13)R}  {}\save[]+<-0.5cm,.5cm>*\txt<8pc>{%
90} \restore & (3,7,9) \ar[r]_*!/_7pt/{(13)P} & (7,9,1) }$$ }
\end{figure}
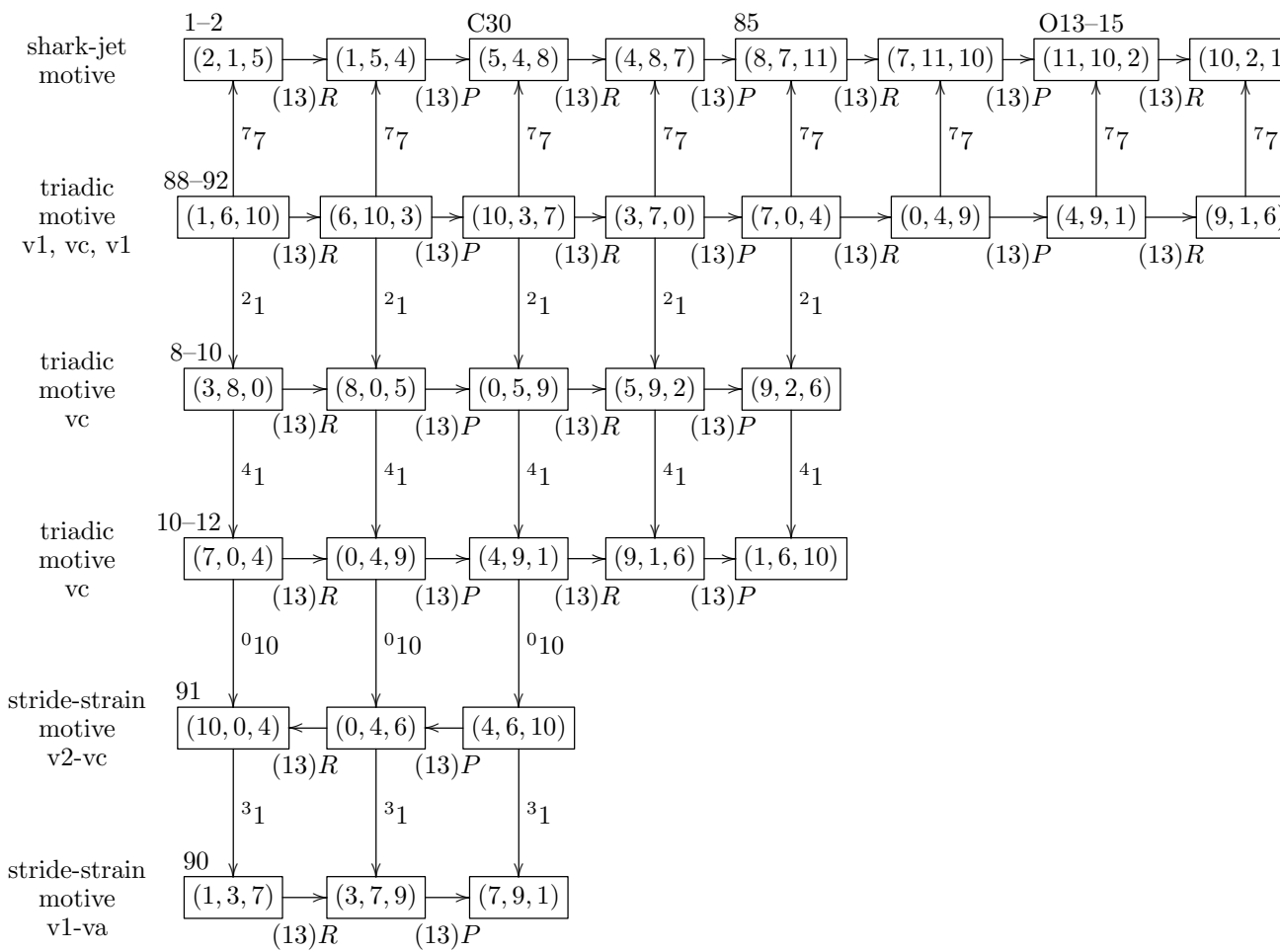

The network in Figure~\ref{fig:Figure8SchoenbergNetworkRP} summarizes the foregoing musical observations. See Section~\ref{sec:Octatonic_Example} for a description of its relationship with Figure~\ref{fig:047_M7_M10_network}. 

It is possible to refine our mathematical framework in order to have all permutations of chords in one GIS, as we do in \cite{fiorenollsatyendraMCM2013}. As an illustration, we now sketch how to include retrogrades  in the Schoenberg example, and explain the Summary Network in Figure~\ref{fig:Figure8SchoenbergNetworkRP}. More precisely, we include the retrograde operation in a pair of dual groups acting on the 48-element union $S$ of the $T/I$-classes of $(7,0,4)$ and $(4,0,7)$. To obtain this pair of dual groups we conjugate the usual $PLR$ operations by two permutations to induce operations on the aforementioned union, and then adjoin the retrograde permutation $(13)$ to the resulting $PLR$-group and the $T/I$-group.

Let $S$ be the set of all transpositions and inversions of both $(7, 0, 4)=(123)(0,4,7)$ and its retrograde $(4, 0, 7)=(12)(0, 4, 7)$.
The notation (123) and (12) {\it without commas} is cycle notation for the permutations $1 \mapsto 2 \mapsto 3 \mapsto 1$ and $1 \mapsto 2 \mapsto 1$. We follow the standard convention for left permutation actions on 3-tuples: $\sigma  (y_1, y_2, y_3 ):=\left(y_{\sigma^{-1}(1)},y_{\sigma^{-1}(2)}, y_{\sigma^{-1}(3)}\right)$. We also follow the usual composition of functions
when we compose permutations, i.e., we do the right function first.
The transpositions and inversions $T_i, I_j \co S \to S$ commute with the permutation action, so in particular they also commute with
the retrograde action $(13)(y_1,y_2,y_3)=(y_3,y_2,y_1)$.
The subgroup of $\text{Sym}(S)$ generated by the $T/I$-group and the permutation $(13)$ is denoted by $\langle T/I, (13) \rangle$. It is
the union
$$\{ T_i \sigma: \sigma =\text{Id}\text{ or } \sigma=(13)\}\cup \{ I_k \sigma:  \sigma =\text{Id}\text{ or } \sigma=(13)\},$$
which clearly has 48 elements. Moreover, it acts simply transitively on $S$. Because $T/I$ and $(13)$ commute, this subgroup $\langle T/I, (13) \rangle\leq \text{Sym}(S)$ is the internal direct product
of the $T/I$-group and the group $\langle (13) \rangle$ in $\text{Sym}(S)$.

To describe the dual group, we define $P$, $L$, and $R$ operations $S \to S$, {\it not uniformly as $J^{1,3}$,  $J^{2,3}$, and $J^{1,2}$,} but rather as
$$(123)P(321), \quad (123)L(321), \quad \text{and} \quad
(123)R(321)$$
on the $T/I$-orbit of $(123)(0,4,7)=(7,0,4)$
and as
$$(12)P(21), \quad (12)L(21), \quad \text{and} \quad
(12)R(21)$$
on the $T/I$-orbit of $(12)(0,4,7)=(4,0,7)$. For simplicity, we just write $P$, $L$, and $R$ for the new operations as well.\footnote{No ambiguity arises because the argument of
the function tells us which conjugate of the usual $P$, $L$, or $R$ to take. }
The other operations in the neo-Riemannian $PLR$-group for $(0,4,7)$ similarly induce operations $S \to S$. Let us temporarily denote this 24-element group by $PLR(S)$. One can check that $PLR(S)$ commutes with the retrograde action of $\{\text{Id},(13)\}$, and also commutes with $T_i$ and $I_j$. Thus, the internal direct product
$$\langle PLR(S), (13) \rangle = \{N\sigma:N\in\text{$PLR(S)$  and } \sigma =\text{Id}\text{ or } \sigma=(13)\}$$
commutes with $\langle T/I,(13) \rangle$, and, as it acts simply transitively on $S$, it is the dual group to $\langle T/I,(13) \rangle$. The intersection of these dual groups is
$$\langle T/I, (13) \rangle\cap\langle PLR(S), (13) \rangle=\langle T_6, (13) \rangle=\{\text{Id}, T_6, (13),T_6(13)\}.$$

We now turn to the second row of the Summary Network. The subgroup of $\langle PLR(S), (13) \rangle$ generated by $(13)R$ and $(13)P$ is dihedral of order 8, and acts simply transitively on its orbit of $(7,0,4)$ by Theorem~\ref{thm:sub_dual_groups}.  This orbit is precisely the second row of the Summary Network. We now have a retrograde refinement of $H_0^{(0,4,7)}$.

The entire previous discussion can be repeated for the jet $(2,1,5)$ and its retrograde $(5,1,2)$ to construct a simply transitive action on their 48 transpositions and inversions. The subgroup generated by $(13)R$ and $(13)P$ is dihedral of order 8, and produces the simply transitive action in the first row of the Summary Network. This subgroup is a retrograde refinement of $H_0^{(1,5,2)}$.

The same goes for the major $(9,2,6)$ and its retrograde $(6,2,9)$ to create row 3 and a refinement of $H_0^{(2,6,9)}$.

Proceeding analogously for the stride $(10,0,4)$ and its retrograde $(4,0,10)$ gives a Klein 4-group generated by $(13)R$ and $(13)P$, rather
than a dihedral group of order 8 (recall the discussion in Section~\ref{sec:Octatonic_Example}). 

Lastly, we need to justify why the affine maps in the Summary Network are morphisms of these refined simply transitive group actions using Theorem~\ref{thm:contextual_inversions_commute_with_affine_maps} and Corollary~\ref{cor:affines_commute_with_subgroups_generated_by_contextual_inversions}, even though the conjugate extensions $P$, $L$, and $R$ are not $J^{1,3}$, $J^{2,3}$, and $J^{1,2}$ globally. Affine maps clearly commute with permutations, so it suffices to prove that they commute with $P$, $L$, and $R$. To exemplify the justification, if $Y$ is a transposed or inverted form of $(0,4,7)$ and $\sigma$ is a permutation, then
$$({}^7 7) P(\sigma Y) \overset{\mathrm{def}}{=}({}^7 7) \sigma J^{1,3}\sigma^{-1}(\sigma Y)=\sigma J^{1,3}\sigma^{-1}(\sigma ({}^7 7)Y)\overset{\mathrm{def}}{=}P(\sigma ({}^7 7) Y),$$
where we have used Theorem~\ref{thm:contextual_inversions_commute_with_affine_maps} to commute ${}^7 7$ past $J^{1,3}$. But this last expression is the same as $P ({}^7 7)\sigma Y$, so that finally $({}^7 7) P= P ({}^7 7)$.


\section{Simply Transitive Covers of the Octatonic}
\label{sec:Covers}

We now explore a kind of converse to the derivation of the octatonic from a complete $PR$-cycle: simply transitive covers of the octatonic.

Beginning with the $PR$-orbit of consonant triads $$\{C,c,E\flat,e\flat,G\flat,g\flat,A,a\},$$ one obtains the octatonic set $O_{01}=\{0,1,3,4,6,7,9,10\}$ as the underlying pitch-class set. This $PR$-orbit forms a {\it cover} of the ocatonic set, and the $PR$-group acts simply transitively on this cover, as does its dual group
\begin{equation} \label{equ:octatonic_setwise_stabilizer}
\{T_0,T_3,T_6,T_9,I_7,I_{10},I_1,I_4\}.
\end{equation} However, instead of starting with an orbit of consonant triads and considering the underlying pitch-class set, the question may be reversed: which (not necessarily consonant) triads in the octatonic set provide a cover with a simply transitive action by the octatonic stabilizer?

Some examples of triads which produce simply transitive covers of the octatonic are listed in Figure~\ref{fig:dual_groups_to_G0}, and diagrammatically depicted in Figure~\ref{fig:047_M7_M10_network}. Via Theorem~\ref{thm:sub_dual_groups}, all three examples in Figure~\ref{fig:dual_groups_to_G0} come from sub dual groups of a duality between a $T/I$-group and a contextual group. However, there are triads in the octatonic which have a $T/I$-symmetry but nevertheless produce a simply transitive cover. In fact, the goal of the current section is to prove that {\it every} triad in the octatonic produces a simply transitive cover, despite any symmetries it may have.

First, we observe that transitivity of the cover follows quickly from transitivity on pitch classes. More precisely, let $O$ be a subset of $\mathbb{Z}_{m}$ (such as the octatonic) and consider the natural mod $m$ $T/I$-action on $\mathbb{Z}_{m}$. Let $G_0$ be any $T/I$-subgroup which preserves $O$ as a set, and let $X$ be any 3-element subset of $O$ (Note: in the previous section $X$ denoted an ordered pitch-class segment, whereas in the present section $X$ denotes an {\it unordered} pitch-class set). If $G_0$ acts transitively on $O$, then $G_0X$ is a cover of $O$. Clearly, $G_0$ acts transitively on $G_0X$. {\it Simple} transitivity of the $G_0$-action requires another hypothesis, as we see in the following observation.

\begin{thm}[Triads generating simply transitive covers] \label{thm:simply_transitive_Zm-covers}
Let $O$ be a subset of $\mathbb{Z}_{m}$ and $X$ any 3-element subset of $O$. Consider the natural mod $m$ $T/I$-action on $\mathbb{Z}_{m}$. Let $G_0$ be any $T/I$-subgroup which preserves $O$ as a set and acts transitively on $O$.

Then $G_0$ acts simply transitively on the transitive $O$-cover $G_0X$ if and only if no nontrivial element of $G_0$ fixes $X$ as a set.
\end{thm}
\begin{pf}
If $G_0$ acts simply transitively, then the stabilizer $(G_0)_X$ is trivial, for $gX=X$ implies $g=e$.

If no nontrivial element of $G_0$ fixes $X$ as a set, and $gX \in G_0X$ and $g_1,g_2 \in G_0$ are such that $g_1(gX)=g_2(gX)$, then $(g_2g)^{-1}(g_1g)$ fixes $X$. We have
$$\aligned
(g_2g)^{-1}(g_1g) &= e \\
g_1g &= g_2g \\
g_1 &= g_2.
\endaligned$$
\end{pf}

In the case of the ocatonic $O_{01}$, the only $T/I$-subgroup which both preserves it and acts transitively is its set-wise stabilizer, displayed in equation \eqref{equ:octatonic_setwise_stabilizer}.
None of its subgroups can act transitively for cardinality reasons.

\begin{cor} \label{cor:simply_transitive_covers_of_octatonic}
Any 3-element subset $X$ of the octatonic $$O_{01}=\{0,1,3,4,6,7,9,10\}\subset \mathbb{Z}_{12}$$ generates a simply transitive cover with respect to the set-wise stabilizer listed in equation \eqref{equ:octatonic_setwise_stabilizer}.
\end{cor}
\begin{pf}
Let $G_0$ denote the set-wise stabilizer of the octatonic $O_{01}$.

It suffices to prove that every 3-element octatonic subset $Y$ with $0 \in Y$ has a trivial $G_0$-stabilizer.
Namely, a general 3-element subset $X$ of the octatonic has a trivial stabilizer if and only if $gX$ does for some (equivalently every) $g \in G_0$ since $(G_0)_{gX}=g(G_0)_{X}g^{-1}$, which means that $X$ generates a simply transitive cover if and only if $gX$ does by Theorem~\ref{thm:simply_transitive_Zm-covers}. The transitivity of $G_0$ on $O_{01}$ guarantees there exists some $g \in G_0$ with $0\in gX$.

Let $Y$ be a 3-element subset of $O_{01}$ with $0 \in Y$. We know $8 \notin Y$ because $8 \notin O_{01}$, so $Y$ is not fixed by $T_3$, $T_6$, or $T_9$ (see Lemma~\ref{lem:048}).

We next need to check if $Y$ is fixed by any of the inversions $I_1$, $I_4$, $I_7$, or $I_{10}$. Recall that if $a,b \in \mathbb{Z}_{12}$, the unique inversion which interchanges $a$ and $b$ is $I_{a+b}$. Our method is to assume $I_j$ interchanges the elements 0 and $b$ of $Y$, which means that $I_j$ fixes the third element $c$, so $c=-c+j$, and then show this leads to an equation in $\mathbb{Z}_{12}$ which cannot be solved in $O_{01}$.  We go through all the pairs $0$ and $b$ which may be in $Y$.

The elements 0 and 0 cannot be interchanged by $I_1$, $I_4$, $I_7$, or $I_{10}$.

The elements 0 and 1 are interchanged by $I_1$, and $c=-c+1$ implies $2c=1$, an equation which cannot be solved in $O_{01}$, and not even in $\mathbb{Z}_{12}$ for parity reasons.

The elements 0 and 3 cannot be interchanged by $I_1$, $I_4$, $I_7$, or $I_{10}$.

The elements 0 and 4 are interchanged by $I_4$, and $c=-c+4$ implies $2c=4$ in $\mathbb{Z}_{12}$, which is $2c-4=12k$ in $\mathbb{Z}$, or equivalently $c-2=6k$ in $\mathbb{Z}$, which implies $c=2$ or $c=8$, neither of which is in $O_{01}$.

The elements 0 and 6 cannot be interchanged by $I_1$, $I_4$, $I_7$, or $I_{10}$.

The elements 0 and 7 are interchanged by $I_7$, and $c=-c+7$ implies $2c=7$, an equation which cannot be solved in $O_{01}$, and not even in $\mathbb{Z}_{12}$ for parity reasons.

The elements 0 and 9 cannot be interchanged by $I_1$, $I_4$, $I_7$, or $I_{10}$.

The elements 0 and 10 are interchanged by $I_{10}$, and $c=-c+10$ implies $2c=10$, which means $c=5$ or 11, neither of which is in $O_{01}$.

In conclusion, $Y$ is not fixed by any nontrivial element of $G_0$, and $Y$ generates a simply transitive cover by Theorem~\ref{thm:simply_transitive_Zm-covers}.
\end{pf}

If the trichord $X$ has no $T/I$-symmetries, then it of course has no $G_0$-symmetries, and $G_0X$ is a simply transitive cover by Theorem~\ref{thm:simply_transitive_Zm-covers}, or by \cite[Theorem~3.1(i)]{fiorenoll2011}. The converse is not true, so Corollary~\ref{cor:simply_transitive_covers_of_octatonic} is new, as the following example illustrates.

\begin{examp}
The chord $\{0,3,6\}$ in the octatonic $O_{01}$ generates a simply transitive cover by Corollary~\ref{cor:simply_transitive_covers_of_octatonic}, although $I_6$ is a symmetry of $\{0,3,6\}$. As above, we denote by $G_0$ the set-wise stabilizer of $O_{01}$ in equation \eqref{equ:octatonic_setwise_stabilizer}. Simple transitivity of $G_0$ on $G_0\{0,3,6\}$ does {\it not} follow from \cite[Theorem~3.1(i)]{fiorenoll2011}, since the $T/I$-group does not act simply transitively on its orbit of $\{0,3,6\}$. We may view the octatonic set $O_{01}$ as a {\it boundary condition} which makes unique solutions in equations for transitivity. For instance, the equation $g\{0,3,6\}=\{3,6,9\}$ does not have a unique solution in the $T/I$-group (e.g. $T_3$ and $T_3I_6$ are both solutions). However, the octatonic set specifies the extra condition that we want a solution which  preserves the octatonic, and this requirement leads to the unique solution $T_3$.
\end{examp}

We finish the proof of Corollary~\ref{cor:simply_transitive_covers_of_octatonic} with the following lemma, which we postponed until now, so as not to distract from the main ideas.

\begin{lem} \label{lem:048}
All 3-element subsets of $\mathbb{Z}_{12}$ which are fixed by some nontrivial translation are of the form $\{k,k+4,k+8\}$ for $k \in \mathbb{Z}_{12}$. Hence, the nontrivial translation must be $T_4$ or $T_8$.
\end{lem}
\begin{pf}
Suppose $a,b,c \in \mathbb{Z}_{12}$ are distinct and $\{a+i,b+i,c+i\}=\{a,b,c\}$. Then $a$ must be equal to one of the elements on the left. If $a=a+i$, then $i=0$ and there is nothing to show.

If $a=b+i$ and $b=a+i$, then $c=c+i$ and $i=0$, again nothing to show.

The other possibility in this second case $a=b+i$ is that $b=c+i$. Then $c=a+i$, and $i$ is the same as
\begin{equation} \label{equ:048}
a-b=b-c=c-a.
\end{equation}
Then $a+c=2b$ and $c=2a-b$. Combining these last two equations, we have $3a=3b$, which is $3(b-a)=0$ in $\mathbb{Z}_{12}$, and therefore $b-a$ is 0, 4, or 8. Similarly, from equation \eqref{equ:048} we see $b+a=2c$, which in combination with $a+c=2b$ gives $3b=3c$, which means $3(c-b)=0$ in $\mathbb{Z}_{12}$, and $c-b$ is 0, 4, or 8. Thus we also have $c-a=(c-b)+(b-a) \in \{0,4,8\}$, and $\{a,b,c\}$ is of the form $\{k,k+4,k+8\}$ for some $k \in \mathbb{Z}_{12}$.

The third possibility for $a$ is that $a=c+i$. But this argument is exactly the same as the two possibilities in the second case $a=b+i$ above, since we may simply relabel the elements $b$ and $c$ to be $c$ and $b$.
\end{pf}

\end{document}